\documentclass[a4paper,12pt]{article} 
\usepackage{amsmath} 
\usepackage{amsopn} 
\usepackage{amssymb} 
\usepackage[mathscr]{eucal} 
\usepackage{theorem} 
\usepackage{enumerate} 

\theorembodyfont{\itshape} 
\theoremstyle{plain}
  \newtheorem{thm}{Theorem}[section] 
  \newtheorem{pro}[thm]{Proposition} 
   
  \newtheorem{cor}[thm]{Corollary} 

\theorembodyfont{\rmfamily} 
\theoremstyle{plain}
  \newtheorem{rem}{Remark}[section] 
  \newtheorem{ex}{Example}[section] 

\renewcommand{\theequation}% 
           {\thesection.\arabic{equation}}

\begin{document} 

\begin{center} 
{\LARGE Sections of time-like twistor spaces with} 

\vspace{1mm} 

{\LARGE light-like or zero covariant derivatives} 

\vspace{6mm} 

{\Large Naoya {\sc Ando}} 
\end{center} 

\vspace{3mm} 

\begin{quote} 
{\footnotesize \it Abstract} \ 
{\footnotesize 
The conformal Gauss maps of time-like minimal surfaces in $E^3_1$ give 
sections of the time-like twistor spaces associated 
with the pull-back bundles 
such that the covariant derivatives are fully light-like, that is, 
these are either light-like or zero, and do not vanish at any point. 
For an oriented neutral $4n$-manifold $(\Tilde{M} , h)$, 
if $J$ is an $h$-reversing almost paracomplex structure of $\Tilde{M}$ 
such that $\nabla J$ is locally given by the tensor product 
of a nowhere zero 1-form and an almost nilpotent structure related to $J$, 
then we will see that $\nabla J$ is valued in 
a light-like $2n$-dimensional distribution $\mathcal{D}$ 
such that $(\Tilde{M} , h, \mathcal{D} )$ is a Walker manifold 
and that the square norm $\parallel\!\nabla J\!\parallel^2$ of $\nabla J$ 
vanishes. 
We will obtain examples of $h$-reversing almost paracomplex structures 
of $E^{4n}_{2n}$ as above. 
In addition, we will obtain 
all the pairs of $h$-reversing almost paracomplex structures of $E^4_2$ 
such that each pair gives sections of the two time-like twistor spaces 
with fully light-like covariant derivatives.} 
\end{quote} 

\vspace{3mm}  

\section{Introduction}\label{sect:intro} 

Neutral metrics are already investigated in various situations. 
The space of oriented lines in the Euclidean $3$-space $E^3$ admits 
a neutral K\"{a}hler structure (\cite{GK}). 
Analogous spaces are found in \cite{AGK}, \cite{salvai}. 
See \cite{salvai2} for fibrations of $E^3$ by oriented lines. 
See \cite{GG} for almost paracomplex structures on neutral 4-manifolds. 
See \cite{dunajski0} for anti-self-dual null-K\"{a}hler structures. 
The ultra-hyperbolic equation is a neutral analogue of the Laplace equation 
(see \cite{asgeirsson}, \cite{CG}, \cite{john}), and 
related to tomography (\cite{guilfoyle}). 
See \cite{pavsic} for quantum field theories in neutral spaces. 

Let $E$ be an oriented vector bundle of rank 4. 
Let $h$ be a positive-definite or neutral metric of $E$ 
and $\nabla$ an $h$-connection of $E$. 
An $h$-preserving complex structure $I$ of $E$ satisfies $\nabla I=0$ 
if and only if 
the corresponding section $\Omega$ of one of 
the twistor or space-like twistor spaces associated with $E$ is 
horizontal with respect to the connection $\hat{\nabla}$ 
induced by $\nabla$. 
If $E$ is the tangent bundle $T\!\Tilde{M}$ of an oriented Riemannian or 
neutral 4-manifold $\Tilde{M}$ and 
if $h$, $\nabla$ are its metric and the Levi-Civita connection of $h$ 
respectively, 
then $I$ is an almost complex structure of $\Tilde{M}$ 
and $\nabla I=0$ just means that $(\Tilde{M} , h, I)$ is 
a K\"{a}hler or neutral K\"{a}hler surface. 
If $E$ is the pull-back bundle $F^*\!T\!\Tilde{M}$ 
over a Riemann surface $M$ 
by a space-like and conformal immersion $F:M\longrightarrow \Tilde{M}$ 
with zero mean curvature vector, 
then $\nabla I=0$ for a twistor lift $\Omega$ of $F$ just means 
that $F$ is isotropic (refer to \cite{friedrich}, \cite{ando4}). 
For the twistor spaces and the space-like twistor spaces, 
refer to \cite{ES}, \cite{BDM2} respectively. 
Suppose that $\Tilde{M}$ is neutral. 
Then, even if the square norm $\parallel\!\nabla I\!\parallel^2$ 
of $\nabla I$ vanishes, $\nabla I$ does not necessarily vanish. 
We say that $(\Tilde{M} , h, I)$ is \textit{isotropic K\"{a}hler\/} 
if $\parallel\!\nabla I\!\parallel^2 =0$. 
See \cite{DDGMMV}, \cite{GM} 
for examples of isotropic K\"{a}hler but non-K\"{a}hler 4-manifolds. 
See \cite{BDM} 
for examples of isotropic K\"{a}hler but non-K\"{a}hler 6-manifolds. 

Let $E$ be as in the beginning of the previous paragraph. 
Let $h$, $\nabla$ be a neutral metric and an $h$-connection of $E$ 
respectively. 
Then an $h$-reversing paracomplex structure $J$ of $E$ 
satisfies $\nabla J=0$ if and only if 
the corresponding section $\Omega$ of one of the time-like twistor spaces 
associated with $E$ is horizontal 
with respect to $\hat{\nabla}$. 
If $E=T\!\Tilde{M}$ for an oriented neutral 4-manifold $\Tilde{M}$, 
then $J$ is an almost paracomplex structure of $\Tilde{M}$ 
and $\nabla J=0$ just means that $(\Tilde{M} , h, J)$ is 
a paraK\"{a}hler surface. 
If $E=F^*\!T\!\Tilde{M}$ 
for a time-like and conformal immersion $F:M\longrightarrow \Tilde{M}$ 
of a Lorentz surface $M$ into $\Tilde{M}$ with zero mean curvature vector 
and if $\Omega$ is a twistor lift of $F$, 
then $\nabla J=0$ implies that $F$ is isotropic. 
Refer to \cite{HM}, \cite{JR} for the time-like twistor spaces. 
Even if $F$ is isotropic, 
$\nabla J$ does not necessarily vanish (\cite{ando4}). 
It is possible that $\hat{\nabla} \Omega$ is either light-like or zero, 
and does not vanish at any point. 
In such a case, 
$\hat{\nabla} \Omega$ determines a light-like one-dimensional subspace 
of the fiber of $\bigwedge^2\!E$ at each point of $M$ 
and then we say that $\hat{\nabla} \Omega$ is \textit{fully light-like}. 
Let $F$ be the conformal Gauss map of a time-like minimal surface 
in $E^3_1$ with nowhere zero curvature. 
Then $\Tilde{M} =S^4_2$, 
and we see that $F$ has zero mean curvature vector and isotropicity 
and that the covariant derivatives of the twistor lifts $\Omega_{\pm}$ 
are fully light-like. 
We will see that for a paracomplex structure $J=J_{\varepsilon}$ 
corresponding to $\Omega_{\varepsilon}$ ($\varepsilon \in \{ +, -\}$), 
$\nabla J$ is locally represented as $\nabla J=\alpha \otimes N$ 
for a nowhere zero $1$-form $\alpha$ and a nilpotent structure $N$ 
related to $J$ and therefore $\nabla J$ is valued in 
a light-like subbundle of the pull-back bundle of rank $2$ 
(Theorem~\ref{thm:cgmtms}). 
We will define a nilpotent structure of 
an oriented vector bundle $E$ of rank $4n$ with a neutral metric $h$. 
Then a nilpotent structure $N$ gives a null structure on each fiber of $E$ 
such that the image is a light-like $2n$-dimensional subspace 
and $h$ is null-Hermitian with respect to $N$ 
(see Section~\ref{sect:nilstr} and refer to \cite{dunajski}). 
We will characterize an $h$-reversing paracomplex structure $J$ of $E$ 
such that for an $h$-connection $\nabla$ of $E$, $\nabla J$ is locally 
represented as $\nabla J=\alpha \otimes N$ for a $1$-form $\alpha$ and 
a nilpotent structure $N$ related to $J$ (Theorem~\ref{thm:hatnabla0g}). 
See \cite{ando7}, \cite{AK} for nilpotent structures for $n=1$. 

Let $\Tilde{M}$ be an oriented neutral $4n$-manifold. 
Then $\Tilde{M}$ has a neutral metric $h$ 
and $h$ gives the Levi-Civita connection $\nabla$. 
Let $J$ be an $h$-reversing almost paracomplex structure of $\Tilde{M}$ 
such that $\nabla J$ is locally represented as $\nabla J=\alpha \otimes N$ 
for a nowhere zero $1$-form $\alpha$ and an almost nilpotent structure $N$ 
related to $J$. 
Then $\nabla J$ is valued in 
a light-like $2n$-dimensional distribution $\mathcal{D}_J$ on $\Tilde{M}$. 
We will see that $(\Tilde{M} , h, \mathcal{D}_J )$ is 
a \textit{Walker manifold}, 
that is, the covariant derivatives of local generators of $\mathcal{D}_J$ 
by any tangent vector of $\Tilde{M}$ are contained in $\mathcal{D}_J$ 
(Theorem~\ref{thm:JDW}). 
This implies that $\mathcal{D}_J$ is involutive 
and therefore $N$ is a null structure. 
See \cite{BGGNV}, \cite{walker} for Walker manifolds 
and see \cite{dunajski} for null structures. 
We will see that the square norm $\parallel\!\nabla J\!\parallel^2$ 
of $\nabla J$ vanishes, that is, $(\Tilde{M} , h, J)$ is 
\textit{isotropic paraK\"{a}hler} (Theorem~\ref{thm:snorm0}). 
See \cite{DDGMMV2} for isotropic paraK\"{a}hler manifolds. 
Suppose $n=1$. 
Let $J$ be an $h$-reversing almost paracomplex structure of $\Tilde{M}$ 
such that the covariant derivative of the corresponding section 
of one of the time-like twistor spaces is fully light-like. 
Then $\nabla J$ is locally represented as above and 
therefore $(\Tilde{M} , h, \mathcal{D}_J )$ is a Walker manifold. 
Let $\mathcal{D}$ be a light-like two-dimensional distribution 
on an oriented neutral $4$-manifold $\Tilde{M}$ 
such that $(\Tilde{M} , h, \mathcal{D} )$ is a Walker manifold. 
Then $\mathcal{D}$ is locally given by ${\rm Im}\,N$ 
for an almost nilpotent structure $N$. 
We will find an $h$-reversing almost paracomplex structure $J$ 
of a neighborhood of each point of $\Tilde{M}$ 
satisfying $\nabla J =\alpha \otimes N$ 
for a nowhere zero $1$-form $\alpha$ and  
that $N$ is related to $J$ (Theorem~\ref{thm:JDW2}). 
In particular, 
for an almost nilpotent structure $N$ of $\Tilde{M}$ 
satisfying $\nabla N=0$, 
there exists $J$ as above on a neighborhood of each point, 
since $(\Tilde{M} , h, \mathcal{D} )$ with $\mathcal{D} ={\rm Im}\,N$ is 
a Walker manifold. 

Let $E$ be an oriented vector bundle 
over $\mbox{\boldmath{$R$}}^m$ ($m\geq 2$) of rank $4n$. 
Let $h$ be a neutral metric of $E$ 
and $\nabla$ a flat $h$-connection of $E$. 
We will find examples of $h$-reversing paracomplex structures of $E$ 
such that for each $J$, 
$\nabla J$ is represented as $\nabla J=\alpha \otimes N$ 
for a nowhere zero $1$-form $\alpha$ and a nilpotent structure $N$ 
related to $J$ (Example~\ref{ex:ex1} $\sim$ Example~\ref{ex:ex4}). 
In particular, we will find examples 
of $h$-reversing almost paracomplex structures of $E^{4n}_{2n}$ 
as above (Remark~\ref{rem:rem}). 
Suppose $n=1$. 
We will characterize sections of the time-like twistor spaces 
associated with $E$ such that the covariant derivatives are 
fully light-like and we will see that  
the covariant derivative of a paracomplex structure $J$ corresponding to 
such a section can be represented as $\nabla J =\alpha \otimes N$ 
for a nowhere zero $1$-form $\alpha$ and a nilpotent structure $N$ 
satisfying $\nabla N=0$ (Proposition~\ref{pro:gsec}). 
In addition, we will obtain all the pairs of sections 
of the two time-like twistor spaces 
with fully light-like covariant derivatives (Theorem~\ref{thm:gsec+-}). 
Then we will find two types of such pairs (Remark~\ref{rem:rem2}). 
One type corresponds to the pair of the lifts of the conformal Gauss map 
of a time-like minimal surface in $E^3_1$ 
and the other type corresponds to the pair of the lifts of 
a time-like surface in a $4$-dimensional neutral space form 
with zero mean curvature vector given in \cite{ando5}. 
Based on these studies, 
we can obtain the results in the case of $E=T\!E^4_2$ 
(Corollary~\ref{cor:E42J}, Remark~\ref{rem:rem2_0}). 
Therefore we can find 
all the pairs of $h$-reversing almost paracomplex structures 
of $E^4_2$ such that each pair gives sections of 
the two time-like twistor spaces 
with fully light-like covariant derivatives. 

\section{Elements of \mbox{\boldmath{$SO(2n, 2n)$}} 
preserving oriented light-like 
\mbox{\boldmath{$2n$}}-planes}\label{sect:so(2n,2n)}  

\setcounter{equation}{0} 

Let $V$ be an oriented $4n$-dimensional vector space 
and $h_V$ a neutral metric of $V$. 
Let $W$ be a light-like $2n$-dimensional subspace of $V$. 
Let $(e_1 , \dots , e_{2n} , e_{2n+1} , \dots , e_{4n} )$ be 
an ordered pseudo-orthonormal basis of $V$ giving the orientation of $V$ 
such that 
$e_1 , \dots , e_{2n}$ (respectively, $e_{2n+1} , \dots , e_{4n}$) 
are space-like (respectively, time-like). 
Suppose that $W$ is spanned by 
\begin{equation} 
\begin{array}{lclcl} 
\xi_1     :=e_1 -e_{2n+1}     , & \ & \dots , & \ & \xi_n    :=e_n -e_{3n} , \\ 
\xi_{n+1} :=e_{n+1} +e_{3n+1} , & \ & \dots , & \ & \xi_{2n} :=e_{2n} +e_{4n} . 
\end{array} 
\label{xi} 
\end{equation} 
Let $T$ be an automorphism of $V$ which preserves $h_V$ and 
the orientation of $V$. 
Then there exists an element $A$ of $SO(2n, 2n)$ satisfying 
\begin{equation*} 
 (Te_1 , \dots , Te_{2n} , Te_{2n+1} , \dots , Te_{4n} )  
=( e_1 , \dots ,  e_{2n} ,  e_{2n+1} , \dots ,  e_{4n} )A.  
\end{equation*} 
We represent $A$ as 
$$A=\left[ 
    \begin{array}{cccc} 
     A_{11} & A_{12} & A_{13} & A_{14} \\ 
     A_{21} & A_{22} & A_{23} & A_{24} \\ 
     A_{31} & A_{32} & A_{33} & A_{34} \\ 
     A_{41} & A_{42} & A_{43} & A_{44} 
      \end{array} 
    \right] ,$$ 
where $A_{ij}$ ($i, j=1, 2, 3, 4$) are $n\times n$ matrices. 
We call $A_{ij}$ the $(i, j)$-\textit{block\/} of $A$. 
We use not only $A_{ij}$ but also $A_{(i, j)}$ 
in order to represent the $(i, j)$-block of $A$. 
We set 
$$A^{\times} 
:=\left[ 
  \begin{array}{cccc} 
   A_{22} & -A_{21} &  A_{42} &  A_{41} \\ 
  -A_{12} &  A_{11} &  A_{32} &  A_{31} \\ 
   A_{24} &  A_{23} &  A_{44} & -A_{43} \\ 
   A_{14} &  A_{13} & -A_{34} &  A_{33} 
    \end{array} 
  \right] .$$ 
Let $T^{\times}$ be an endomorphism of $V$ 
such that $A^{\times}$ is the representation matrix with respect to 
$(e_1 , \dots , e_{2n} , e_{2n+1} , \dots , e_{4n} )$. 

\begin{pro}\label{pro:2ndim} 
The following are mutually equivalent\/$:$ 
\begin{itemize} 
\item[{\rm (a)}]{$T$ induces an endomorphism $\Phi_T$ of $W;$} 
\item[{\rm (b)}]{$A_{(i, j)}   +(-1)^j A_{(i,   j+2)} 
         =(-1)^i (A_{(i+2, j)} +(-1)^j A_{(i+2, j+2)})$ 
for $i, j=1, 2;$} 
\item[{\rm (c)}]{$A_{(i, j)}   +(-1)^{i-1} A_{(i+2, j)} 
     =(-1)^{j-1} (A_{(i, j+2)} +(-1)^{i-1} A_{(i+2, j+2)})$ 
for $i, j=1, 2;$} 
\item[{\rm (d)}]{$T^{\times}$ induces 
an endomorphism $\Phi_{T^{\times}}$ of $W$.} 
\end{itemize} 
\end{pro} 

\vspace{3mm} 

\par\noindent 
\textit{Proof} \ 
We see that (a) and (b) are equivalent. 
We can rewrite (b) into (c) immediately. 
We see that (c) and (d) are equivalent. 
\hfill 
$\square$ 

\vspace{3mm} 

For $i, j=1, 2$, we set 
\begin{equation*}
\begin{split} 
& P_{ij}          :=A_{(i, j)}   +(-1)^j A_{(i,   j+2)} , \\  
& P^{\times}_{ij} :=(-1)^{i'+j'} 
                    (A_{(i', j')}   +(-1)^{i'-1} A_{(i'+2, j')} ),  
\end{split} 
\end{equation*} 
where $\{ i, i'\} =\{ j, j'\} =\{1, 2\}$. 
We set 
$$P:=\left[ \begin{array}{cc} 
             P_{11} & P_{12} \\ 
             P_{21} & P_{22} 
              \end{array} 
     \right] , \quad 
  P^{\times} 
:=\left[ \begin{array}{cc} 
          P^{\times}_{11} & P^{\times}_{12} \\ 
          P^{\times}_{21} & P^{\times}_{22} 
           \end{array} 
  \right] .$$ 
If $T$ induces an endomorphism $\Phi_T$ of $W$, 
that is, 
if $T^{\times}$ induces an endomorphism $\Phi_{T^{\times}}$ of $W$, 
then $P$, $P^{\times}$ are the representation matrices 
of $\Phi_T$, $\Phi_{T^{\times}}$ 
respectively with respect to $\xi_1 , \dots , \xi_{2n}$. 
Let $\Lambda_n$ be a $4n\times 4n$ matrix defined by 
$$\Lambda_n 
:=\left[ \begin{array}{cccc} 
          O_n & -I_n & O_n &  I_n \\ 
          I_n &  O_n & I_n &  O_n \\ 
          O_n &  I_n & O_n & -I_n \\ 
          I_n &  O_n & I_n & O_n 
           \end{array} 
  \right] ,$$ 
where $I_n$ is the $n\times n$ unit matrix 
and   $O_n$ is the $n\times n$ zero matrix. 
We obtain 

\begin{pro}\label{pro:commute} 
An automorphism $T$ satisfies not only 
one of {\rm (a)}$\sim${\rm (d)} in Proposition~\ref{pro:2ndim} 
but also $P=P^{\times}$ 
if and only if $A$ and $\Lambda_n$ are commutative. 
\end{pro} 

\begin{rem} 
The condition $A\Lambda_n =\Lambda_n A$ just means 
that $T$ is null-linear (see \cite{dunajski}). 
\end{rem} 

Let $G$ be a Lie subgroup of $SO(2n, 2n)$ 
defined by $A\Lambda_n =\Lambda_n A$ for $A\in SO(2n, 2n)$. 
Then each element of the Lie algebra of $G$ is represented as 
$$\left[ 
  \begin{array}{cccc} 
   U_{11}            &   U_{12}    &   
   U_{11} +Z         &  -U_{12} +X       \\ 
  -{}^t U_{12}       &   U_{22}    &  
  -{}^t U_{12} +Y    &  -U_{22} +{}^t Z \\ 
  -U_{11} +{}^t Z    &  -U_{12} +Y &  
  -U_{11} -Z +{}^t Z &   U_{12} -X -Y    \\      
  -{}^t U_{12} +X    &   U_{22} +Z &  
  -{}^t U_{12} +X +Y &  -U_{22} -Z +{}^t Z 
    \end{array}  
  \right] ,$$ 
where 
$${}^t  U_{11} =-U_{11} , \quad 
  {}^t  U_{22} =-U_{22} , \quad  
  {}^t\!X      =X,        \quad 
  {}^t  Y      =Y.$$  
Therefore we obtain $\dim G=4n^2$. 
We will prove 

\begin{pro}\label{pro:2ndim2} 
Suppose that $T$ satisfies 
one of {\rm (a)}$\sim${\rm (d)} in Proposition~\ref{pro:2ndim}. 
Then 
\begin{itemize} 
\item[{\rm (a)}]{$|P||P^{\times}|=1;$} 
\item[{\rm (b)}]{if $P=P^{\times}$, that is, if $A\in G$, 
then $|P|=1$, that is, $\xi :=\xi_1 \wedge \dots \wedge \xi_{2n}$ is 
invariant by $\Phi_T =\Phi_{T^{\times}}$.} 
\end{itemize} 
\end{pro} 

\vspace{3mm} 

\par\noindent 
\textit{Proof} \ 
Since $A\in SO(2n, 2n)$, we have $|A|=1$. 
Using (b) in Proposition~\ref{pro:2ndim}, we obtain 
\begin{equation*} 
\begin{split} 
& 
\left| \begin{array}{cccc} 
        A_{11} & A_{12} & A_{13} & A_{14} \\ 
        A_{21} & A_{22} & A_{23} & A_{24} \\ 
        A_{31} & A_{32} & A_{33} & A_{34} \\ 
        A_{41} & A_{42} & A_{43} & A_{44} 
         \end{array} 
\right| \\ 
& = 
\left| \begin{array}{cccc} 
        A_{11} &  A_{12} &  A_{13} -A_{11} &  A_{14} +A_{12} \\ 
        A_{21} &  A_{22} &  A_{23} -A_{21} &  A_{24} +A_{22} \\ 
        A_{31} &  A_{32} &  A_{33} -A_{31} &  A_{34} +A_{32} \\ 
        A_{41} &  A_{42} &  A_{43} -A_{41} &  A_{44} +A_{42} 
         \end{array} 
\right| \\ 
& = 
\left| \begin{array}{cccc} 
        A_{11} &  A_{12} &  A_{13} -A_{11} &   A_{14} +A_{12} \\ 
        A_{21} &  A_{22} &  A_{23} -A_{21} &   A_{24} +A_{22} \\ 
        A_{31} &  A_{32} &  A_{11} -A_{13} &  -A_{14} -A_{12} \\ 
        A_{41} &  A_{42} &  A_{23} -A_{21} &   A_{22} +A_{24} 
         \end{array} 
\right| \\ 
& = 
\left| \begin{array}{cccc} 
        A_{11} +A_{31} &  A_{12} +A_{32} &  O_n            &   O_n \\ 
        A_{21} -A_{41} &  A_{22} -A_{42} &  O_n            &   O_n \\ 
        A_{31}         &  A_{32}         &  A_{11} -A_{13} &  -A_{14} -A_{12} 
        \\ 
        A_{41}         &  A_{42}         &  A_{23} -A_{21} &   A_{22} +A_{24} 
         \end{array} 
\right| . 
\end{split} 
\end{equation*} 
Therefore we obtain (a) in Proposition~\ref{pro:2ndim2}. 
We will prove (b) in Proposition~\ref{pro:2ndim2} 
by induction as follows. 
\\ 
\textit{Part 1} \ 
Suppose $n=1$. 
Then an element $A$ of $SO(2, 2)$ is represented as $A=BC$, 
where 
$$B=\left[ \begin{array}{cccc} 
            b_1 & -b_2 &  b_3 &  b_4 \\ 
            b_2 &  b_1 & -b_4 &  b_3 \\ 
            b_3 & -b_4 &  b_1 &  b_2 \\ 
            b_4 &  b_3 & -b_2 &  b_1 
             \end{array} 
           \right]  \quad 
\left( 
\begin{array}{l} 
  b_1 ,  b_2 ,  b_3 ,  b_4 \in \mbox{\boldmath{$R$}}, \\ 
  b^2_1 +b^2_2 -b^2_3 -b^2_4 =1 
  \end{array} 
\right)$$ 
and 
$$C=\left[ \begin{array}{cccc} 
            1 & 0      & 0      & 0      \\ 
            0 & c_{22} & c_{23} & c_{24} \\ 
            0 & c_{32} & c_{33} & c_{34} \\ 
            0 & c_{42} & c_{43} & c_{44} 
             \end{array} 
    \right] \quad 
([c_{ij} ]\in SO(1, 2)).$$ 
Suppose $C=I_4$. Then $A=B$ satisfies (b) in Proposition~\ref{pro:2ndim}, 
$P=P^{\times}$ and $|P|=1$. 
Suppose $B=I_4$. Then $A=C$ satisfies (b) in Proposition~\ref{pro:2ndim} 
if and only if $[c_{ij} ]\in SO(1, 2)$ satisfies 
\begin{equation} 
\left[ \begin{array}{ccc} 
        c_{22} & c_{23} & c_{24} \\ 
        c_{32} & c_{33} & c_{34} \\ 
        c_{42} & c_{43} & c_{44} 
         \end{array} 
\right]  
= 
\left[ \begin{array}{ccc} 
 \varepsilon ({\rm cosh}\,t +e^{-t} c^2 /2) &  c &   
 \varepsilon ({\rm sinh}\,t -e^{-t} c^2 /2) \\ 
 \varepsilon e^{-t} c                       &  1 &  
-\varepsilon e^{-t} c                       \\ 
 \varepsilon ({\rm sinh}\,t +e^{-t} c^2 /2) &  c &  
 \varepsilon ({\rm cosh}\,t -e^{-t} c^2 /2)  
         \end{array} 
\right] , 
\label{c22244244} 
\end{equation} 
where $\varepsilon =+1$ or $-1$, and $c$, $t\in \mbox{\boldmath{$R$}}$. 
Suppose that $[c_{ij} ]\in SO(1, 2)$ satisfies \eqref{c22244244}. 
Then we have 
$$P=\left[ 
    \begin{array}{cc} 
     1 &   0 \\ 
    -c &  \varepsilon e^t 
      \end{array} 
    \right] , 
    \quad 
  P^{\times} 
   =\left[ 
    \begin{array}{cc} 
     \varepsilon e^{-t}   &  0 \\ 
    -\varepsilon e^{-t} c &  1 
      \end{array} 
    \right] .$$ 
Therefore the following are mutually equivalent: 
(i)   $\varepsilon =1$ and $t=0$; 
(ii)  $P=P^{\times}$; 
(iii) $|P|=|P^{\times}|=1$. 
In particular, if $P=P^{\times}$, then $|P|=1$. 
In general, 
if $A=BC$ satisfies not only 
one of {\rm (a)}$\sim${\rm (d)} in Proposition~\ref{pro:2ndim} 
but also $P=P^{\times}$, 
then noticing Proposition~\ref{pro:commute}, 
we see that both $B$ and $C$ satisfy these conditions. 
Therefore we obtain $|P|=1$ for $T$ given by $A=BC$. 
Hence we obtain (b) in Proposition~\ref{pro:2ndim2} 
in the case of $n=1$. \\ 
\textit{Part 2} \ 
Suppose $n\geq 2$. 
We set 
\begin{equation} 
H:=\left\{ \left[ 
           \begin{array}{cccc} 
            A_{11} & -A_{21} &  A_{31} & A_{41} \\ 
            A_{21} &  A_{11} & -A_{41} & A_{31} \\ 
            A_{31} & -A_{41} &  A_{11} & A_{21} \\ 
            A_{41} &  A_{31} & -A_{21} & A_{11} 
             \end{array} 
           \right] \in SO(2n, 2n) 
   \right\} .  
\label{H} 
\end{equation} 
Then $H$ is a Lie subgroup of $G$ with $\dim H =2n^2 +n$. 
Since there exists a path in $H$ 
from the unit element to each element of $H$, 
any element of $H$ satisfies $|P|=1$. 
For each element $A$ of $G$, 
there exists an element $B$ of $H$ 
such that $C:=AB$ satisfies  
\begin{equation*} 
\begin{split} 
C_{11} 
=\left[ \begin{array}{cccc} 
          1     &   0     &  \cdots &   0     \\ 
          0     &   *     &  \cdots &   *     \\ 
         \vdots &  \vdots &  \ddots &  \vdots \\ 
          0     &   *     &  \cdots &   * 
          \end{array} 
 \right] , \ 
C_{1j} & 
=\left[ \begin{array}{cccc} 
          0     &   0     &  \cdots &   0     \\ 
          *     &   *     &  \cdots &   *     \\ 
         \vdots &  \vdots &  \ddots &  \vdots \\ 
          *     &   *     &  \cdots &   * 
          \end{array} 
 \right] \ (j=2, 3, 4), \\ 
C_{i1} & 
=\left[ \begin{array}{cccc} 
          0     &   *     &  \cdots &   *     \\ 
          0     &   *     &  \cdots &   *     \\ 
         \vdots &  \vdots &  \ddots &  \vdots \\ 
          0     &   *     &  \cdots &   * 
           \end{array} 
  \right] \ (i=2, 3, 4). 
\end{split} 
\end{equation*} 
In addition, noticing $A\Lambda_n =\Lambda_n A$ and $A\in SO(2n, 2n)$, 
we see that $C$ satisfies 
\begin{equation*} 
C_{ii}  
=\left[ \begin{array}{cccc} 
         1     &   0     &  \cdots &   0     \\ 
         0     &   *     &  \cdots &   *     \\ 
        \vdots &  \vdots &  \ddots &  \vdots \\ 
         0     &   *     &  \cdots &   * 
          \end{array} 
 \right] \ (i=2, 3, 4) 
\end{equation*} 
and 
\begin{equation*} 
C_{ij} 
=\left[ \begin{array}{cccc} 
         0     &   0     &  \cdots &   0     \\ 
         0     &   *     &  \cdots &   *     \\ 
        \vdots &  \vdots &  \ddots &  \vdots \\ 
         0     &   *     &  \cdots &   * 
          \end{array} 
 \right] \ (i\not= j). 
\end{equation*} 
Therefore noticing Part 1, 
we see by induction that $A$ satisfies $|P|=1$. 
Hence we obtain (b) in Proposition~\ref{pro:2ndim2} 
in the case of $n\geq 2$. 
\hfill 
$\square$ 

\vspace{3mm} 

\begin{rem}\label{rem:|P|=1}
Suppose $n\geq 2$ and that $A\in SO(2n, 2n)$ satisfies 
one of {\rm (a)}$\sim${\rm (d)} in Proposition~\ref{pro:2ndim}. 
Then even if $|P|=1$, it is possible that $P\not= P^{\times}$. 
For example, 
if $X$, $Y$ are distinct elements of $SO(n)$, 
then 
$$\left[ \begin{array}{cccc} 
          X   & O_n & O_n & O_n \\ 
          O_n & Y   & O_n & O_n \\ 
          O_n & O_n & X   & O_n \\ 
          O_n & O_n & O_n & Y 
           \end{array} 
  \right]$$ 
satisfies one of {\rm (a)}$\sim${\rm (d)} in Proposition~\ref{pro:2ndim} 
and $|P|=1$ but does not satisfy $P\not= P^{\times}$. 
\end{rem} 

\begin{rem}\label{rem:surj} 
We set 
$$\Tilde{G} 
:=\{ A\in GL(4n, \mbox{\boldmath{$R$}}) \ | \ 
     A\Lambda_n =\Lambda_n A\} .$$  
Then $\Tilde{G}$ is a Lie subgroup of $GL(4n, \mbox{\boldmath{$R$}})$ 
with $\dim \Tilde{G} =8n^2$ 
and $G$ is a Lie subgroup of $\Tilde{G}$. 
Referring to the proof of (a) in Proposition~\ref{pro:2ndim2}, 
we see that $\Tilde{G}$ is contained in 
the connected component $GL_0 (4n, \mbox{\boldmath{$R$}})$ 
of $GL(4n, \mbox{\boldmath{$R$}})$ with the unit element. 
Since each $A\in \Tilde{G}$ gives an automorphism of $W$, 
we obtain a homomorphism $\psi$ from $\Tilde{G}$ 
to $GL(2n, \mbox{\boldmath{$R$}})$. 
Since $\dim {\rm Ker}\,\psi =4n^2$, 
$\psi$ is a surjective homomorphism 
from the connected component of $\Tilde{G}$ with the unit element 
onto $GL_0 (2n, \mbox{\boldmath{$R$}})$. 
In addition, 
if $A$ is an element of $\Tilde{G}$ given by 
$$A=\left[ \begin{array}{cccc} 
            I_n & O_n & O_n      &   O_n      \\ 
            O_n & I_n & O_n      & -2Z_n      \\ 
          -2Z_n & O_n & I_{n, -} &   O_n      \\ 
            O_n & O_n & O_n      &   I_{n, -}   
             \end{array} 
    \right] ,$$ 
where 
\begin{equation} 
I_{1, -} :=-1, \quad 
I_{n, -} :=\left[ \begin{array}{cccc} 
                   -1     &   0     &  \cdots &   0     \\ 
                    0     &   1     &  \ddots &  \vdots \\ 
                   \vdots &  \ddots &  \ddots &   0     \\ 
                    0     &  \cdots &   0     &   1   
                    \end{array} 
           \right] \ (n\geq 2)  
\label{IIZ}
\end{equation}
and 
$$Z_1 := 1, \quad 
  Z_n :=\left[ \begin{array}{cccc} 
                1     &   0     &  \cdots &   0     \\ 
                0     &   0     &  \cdots &   0     \\ 
               \vdots &  \vdots &  \ddots &  \vdots \\ 
                0     &   0     &  \cdots &   0    
                 \end{array} 
        \right] \ (n\geq 2),$$  
then $P$ satisfies $P_{11} =I_n$, $P_{22} =I_{n, -}$ and 
$P_{12} =P_{21} =O_n$. 
Therefore $\psi : \Tilde{G} \longrightarrow GL(2n, \mbox{\boldmath{$R$}})$ 
is surjective. 
We see from (b) in Proposition~\ref{pro:2ndim2} that 
${\rm Im}\,\psi |_G$ is contained in $SL(2n, \mbox{\boldmath{$R$}})$, 
and since $\dim {\rm Ker}\,\psi |_G =2n^2 -n$, 
$\psi |_G : G\longrightarrow SL(2n, \mbox{\boldmath{$R$}})$ is not surjective 
if $n\geq 2$. 
\end{rem} 

\section{Nilpotent structures}\label{sect:nilstr}  

\setcounter{equation}{0} 

Let $M$ be a manifold 
and $E$ an oriented vector bundle over $M$ of rank $4n$. 
Let $h$ be a neutral metric of $E$. 
Let $N$ be a section of ${\rm End}\,E$. 
We say that $N$ is a \textit{nilpotent structure\/} of $E$ 
if on a neighborhood of each point of $M$, 
there exists an ordered pseudo-orthonormal local frame field 
$e=(e_1 , \dots , e_{2n} , e_{2n+1} , \dots , e_{4n} )$ of $E$ 
satisfying $Ne=e\Lambda_n$ and 
that $e_1 , \dots , e_{2n}$ (respectively, $e_{2n+1} , \dots , e_{4n}$) 
are space-like (respectively, time-like). 
Let $N$ be a nilpotent structure of $E$. 
Then we have 
\begin{itemize} 
\item[{\rm (i)}]{${\rm Im}\,N={\rm Ker}\,N$,} 
\item[{\rm (ii)}]{$\pi_N :={\rm Im}\,N={\rm Ker}\,N$ is a subbundle of $E$ 
of rank $2n$ such that each fiber is light-like,} 
\item[{\rm (iii)}]{$h(\phi , N\phi )=0$ for any local section $\phi$ of $E$.} 
\end{itemize} 
In particular, $N$ gives a null structure on each fiber of $E$ 
and $h$ is null-Hermitian with respect to $N$ (see \cite{dunajski}). 

Let $N$ be a section of ${\rm End}\,E$ satisfying (i), (ii), (iii) 
in the previous paragraph. 
Then there exist local sections $\xi_1 , \dots , \xi_{2n}$ of $\pi_N$ 
which form a local frame field of $\pi_N$ 
and local sections $\xi'_1 , \dots , \xi'_{2n}$ of $E$ satisfying 
\begin{itemize} 
\item[{\rm (i)}]{$h(\xi'_i , \xi'_j )=0$ for $i, j=1, \dots , 2n$,} 
\item[{\rm (ii)}]{$h(\xi_i , \xi'_j )=0$ for $i, j=1, \dots , 2n$ 
with $i\not= j$,} 
\item[{\rm (iii)}]{$N(\xi'_i )=\xi_{n+i}$, $N(\xi'_{n+i} )=\xi_i$ 
for $i=1, \dots , n$,} 
\item[{\rm (iv)}]{for $i=1, \dots , n$, 
$h(\xi_i , \xi'_i )=1$, that is, $h(\xi_{n+i} , \xi'_{n+i} )=-1$.} 
\end{itemize} 
For $i=1, \dots , n$, we set 
\begin{equation*} 
\begin{array}{lcl} 
  e_i      :=\dfrac{1}{2} ( 2\xi'_i     +\xi_i ), & \ &  
  e_{n+i}  :=\dfrac{1}{2} (-2\xi'_{n+i} +\xi_{n+i} ), \\ 
  e_{2n+i} :=\dfrac{1}{2} ( 2\xi'_i     -\xi_i ), & \ &  
  e_{3n+i} :=\dfrac{1}{2} ( 2\xi'_{n+i} +\xi_{n+i} ). 
  \end{array} 
\end{equation*} 
Then $(e_1 , \dots , e_{2n} , e_{2n+1} , \dots , e_{4n} )$ is 
an ordered pseudo-orthonormal local frame field of $E$ 
such that 
$e_1 , \dots , e_{2n}$ (respectively, $e_{2n+1} , \dots , e_{4n}$) are 
space-like (respectively, time-like). 
From $N(\xi'_i )=\xi_{n+i}$, $N(\xi'_{n+i} )=\xi_i$, 
we see that $N$ is a nilpotent structure of $E$. 
Let $S$ be an element of $GL(4n, \mbox{\boldmath{$R$}} )$ satisfying 
$$(e_1 , \dots , e_{4n} ) 
 =(\xi'_1 , \dots , \xi'_{2n} , \xi_1 , \dots , \xi_{2n} )S.$$ 
Then $\det S>0$. 
Therefore $(\xi'_1 , \dots , \xi'_{2n} , \xi_1 , \dots , \xi_{2n} )$ gives 
the orientation of $E$ 
if and only if $(e_1 , \dots , e_{4n} )$ gives the orientation of $E$. 

Let $N$ be a nilpotent structure of $E$. 
For $\varepsilon \in \{ +, -\}$, we set 
$$I'_{4n, \varepsilon} 
:=\left[ \begin{array}{cccc} 
          I_n & O_n & O_n & O_n \\ 
          O_n & I_n & O_n & O_n \\ 
          O_n & O_n & I_n & O_n \\ 
          O_n & O_n & O_n & I_{n, \varepsilon} 
           \end{array} 
  \right] ,$$
where $I_{n, +} :=I_n$ and $I_{n, -}$ is as in \eqref{IIZ}.  
For $\varepsilon \in \{ +, -\}$, 
we say that $N$ is an $\varepsilon$-\textit{nilpotent structure\/}  
if we can choose an ordered pseudo-orthonormal local frame field $e$ 
of $E$ on a neighborhood of each point of $M$ 
giving the orientation of $E$ and 
satisfying $NeI'_{4n, \varepsilon} =eI'_{4n, \varepsilon} \Lambda_n$. 
Let $N$ be an $\varepsilon$-nilpotent structure of $E$. 
Then such a frame field as $e$ for $N$ is called 
an \textit{admissible frame field\/} of $N$. 
For an admissible frame field $e$ of $N$, 
\begin{equation} 
  \begin{array}{lcl} 
   \xi_1      :=e_1     -            e_{2n+1} , & \ &  
   \xi_i      :=e_i     -            e_{2n+i} , \\ 
   \xi_{n+1}  :=e_{n+1} +\varepsilon e_{3n+1} , & \ &  
   \xi_{n+i}  :=e_{n+i} +            e_{3n+i}  
    \end{array} \ 
(i=2, \dots , n) 
    \label{xi2} 
\end{equation} 
form a local frame field of $\pi_N$ (\eqref{xi2} with $\varepsilon =+$ 
coincides with \eqref{xi}). 
Let $f=(f_1 , \dots , f_{2n} , f_{2n+1} , \dots , f_{4n} )$ be 
an ordered pseudo-orthonormal local frame field of $E$ 
giving the orientation. 
Then $f$ is an admissible frame field of $N$ 
if and only if for each admissible frame field $e$ of $N$, 
an $SO(2n, 2n)$-valued function $A$ 
on the intersection of the domains of $e$ and $f$ 
given by $fI'_{4n, \varepsilon} =eI'_{4n, \varepsilon} A$ is valued in $G$.  
Therefore we see by (b) in Proposition~\ref{pro:2ndim2} 
that for $\xi_1 , \dots , \xi_{2n}$ as in \eqref{xi2}, 
$\xi =\xi_1 \wedge \dots \wedge \xi_{2n}$ does not depend on the choice 
of an admissible frame field of $N$ 
and this means that $N$ gives a section $\xi$ of $\bigwedge^{2n}\!E$. 

Let $\xi$ be a section of $\bigwedge^{2n}\!E$ 
and suppose that $\xi$ is locally represented 
as $\xi =\xi_1 \wedge \dots \wedge \xi_{2n}$ 
for local light-like sections $\xi_1 , \dots , \xi_{2n}$ of $E$ 
defined on a neighborhood $U$ of each point of $M$ 
such that at each point of $U$, 
they span a light-like $2n$-dimensional subspace of the fiber of $E$. 
Then $\xi_1 , \dots , \xi_{2n}$ are represented 
as in \eqref{xi2} for $\varepsilon \in \{ +, -\}$ and 
an ordered pseudo-orthonormal local frame field 
$e=(e_1 , \dots , e_{2n} , e_{2n+1} , \dots , e_{4n} )$ of $E$ 
giving the orientation of $E$. 
Suppose that $n$ is even. 
Then noticing that there exists $A\in SO(2n, 2n)$ satisfying 
one of {\rm (a)}$\sim${\rm (d)} in Proposition~\ref{pro:2ndim}, 
$|P|=1$ and $P\not= P^{\times}$ (Remark~\ref{rem:|P|=1}), 
we see that $\xi$ gives plural $\varepsilon$-nilpotent structures 
of the restriction of $E$ on a neighborhood of each point of $M$. 
Suppose that $n$ is odd. 
Then $\xi$ gives at least one $\varepsilon$-nilpotent structure 
as above and we see in particular that $n=1$ if and only if 
$\xi$ determines a unique $\varepsilon$-nilpotent structure of $E$. 
See \cite{ando7}, \cite{AK} for the case of $n=1$. 

In the rest of this section, 
suppose $n=1$. 
Let $\hat{h}$ be the metric of $\bigwedge^2\!E$ induced by $h$. 
Then $\hat{h}$ has signature $(2,4)$. 
Noticing 
the double covering $SO_0(2,2)\longrightarrow SO_0 (1,2)\times SO_0(1,2)$, 
we see that $\bigwedge^2\!E$ is represented as the direct sum of 
its two orientable subbundles $\bigwedge^2_+\!E$, $\bigwedge^2_-\!E$ 
of rank $3$. 
We see that $\bigwedge^2_+\!E$ is orthogonal to $\bigwedge^2_-\!E$ 
with respect to $\hat{h}$ 
and that the restriction of $\hat{h}$ on 
each of $\bigwedge^2_+\!E$, $\bigwedge^2_-\!E$ has signature $(1,2)$. 
Let $(e_1 , e_2 , e_3 , e_4 )$ be 
an ordered pseudo-orthonormal local frame field of $E$ 
giving the orientation of $E$. 
We set 
\begin{equation} 
\begin{split} 
\Omega_{\pm , 1} 
& :=\dfrac{1}{\sqrt{2}} (e_1 \wedge e_2 \pm e_3 \wedge e_4 ), \\ 
\Omega_{\pm , 2} 
& :=\dfrac{1}{\sqrt{2}} (e_1 \wedge e_3 \pm e_4 \wedge e_2 ), \\ 
\Omega_{\pm , 3} 
& :=\dfrac{1}{\sqrt{2}} (e_1 \wedge e_4 \pm e_2 \wedge e_3 ). 
\end{split} 
\label{Omega+-123} 
\end{equation}
Then $\Omega_{\pm , 1}$ are space-like 
and  $\Omega_{\pm , 2}$, $\Omega_{\pm , 3}$ are time-like, 
and we can suppose that $\Omega_{-, 1}$, $\Omega_{+, 2}$, $\Omega_{+, 3}$ 
(respectively,          $\Omega_{+, 1}$, $\Omega_{-, 2}$, $\Omega_{-, 3}$) 
form a pseudo-orthonormal local frame field of $\bigwedge^2_+\!E$ 
(respectively, $\bigwedge^2_-\!E$). 
The \textit{light-like twistor spaces\/} associated with $E$ are 
fiber bundles 
$$U_0\!\left(\textstyle\bigwedge^2_{\pm}\!E\right) 
     :=\left\{ \left. \theta \in \textstyle\bigwedge^2_{\pm}\!E 
                                 \setminus \{ 0\} 
     \ \right| \ \hat{h} (\theta , \theta )=0\right\}$$ 
in $\bigwedge^2_{\pm}\!E$ respectively. 
We can refer to \cite{ando7}, \cite{AK} for the light-like twistor spaces. 
Let $\Omega$ be a section of $U_0\!\left(\bigwedge^2_{\varepsilon}\!E
                                   \right)$ ($\varepsilon \in \{ +, -\}$). 
Then we can find $(e_1 , e_2 , e_3 , e_4 )$ as above satisfying 
\begin{equation*} 
\Omega =\Omega_{-\varepsilon , 1} +\varepsilon \Omega_{\varepsilon , 3} . 
\label{Omega+-} 
\end{equation*} 
Therefore $\Omega$ is locally represented 
as $\Omega =(1/\sqrt{2} )\xi_1 \wedge \xi_2$. 
Then $\Omega$ gives a unique $\varepsilon$-nilpotent structure $N$ of $E$. 
It satisfies 
\begin{equation} 
  (Ne_1 , Ne_2 , Ne_3 , N\varepsilon e_4 ) 
 =( e_1 ,  e_2 ,  e_3 ,  \varepsilon e_4 )\Lambda_1 , 
\label{N} 
\end{equation} 
that is, $(e_1 , e_2 , e_3 , e_4 )$ is an admissible frame field of $N$. 
Notice that in \eqref{N}, $\varepsilon$ is put before $e_4$, 
while in the corresponding equation in \cite{ando7}, 
$\varepsilon$ is put before $e_3$.  
Each $\varepsilon$-nilpotent structure of $E$ gives a unique section 
of $U_0\!\left(\bigwedge^2_{\varepsilon}\!E\right)$ and therefore 
there exists a one-to-one correspondence between the set of 
sections of $U_0\!\left(\bigwedge^2_{\varepsilon}\!E\right)$ and 
the set of $\varepsilon$-nilpotent structures of $E$. 

\section{The covariant derivatives of 
paracomplex structures}\label{sect:cdpsls2} 

\setcounter{equation}{0} 

Let $E$, $h$ be as in the beginning of the previous section. 
Let $\nabla$ be an $h$-connection of $E$, i.e., 
a connection of $E$ satisfying $\nabla h=0$. 
Let $J$ be a section of ${\rm End}\,E$. 
For $\varepsilon \in \{ +, -\}$, 
we say that $J$ is an $\varepsilon$-\textit{paracomplex structure\/} 
of $E$ if $J$ satisfies 
\begin{itemize} 
\item[{\rm (i)}]{$J$ is a paracomplex structure of $E$,}
\item[{\rm (ii)}]{$J$ is $h$-reversing, that is, $J^* h=-h$,} 
\item[{\rm (iii)}]{on a neighborhood of each point of $M$, 
there exists an ordered pseudo-orthonormal local frame field 
$e=(e_1 ,\dots e_{2n} , e_{2n+1} , \dots , e_{4n} )$ of $E$ 
giving the orientation and satisfying 
\begin{equation} 
\begin{split}
& (Je_1 , \dots , Je_{2n} , Je_{2n+1} , \dots , Je_{4n} ) 
   I'_{4n, \varepsilon} \\ 
& =(e_1 , \dots ,  e_{2n} ,  e_{2n+1} , \dots ,  e_{4n} ) 
   I'_{4n, \varepsilon}  
  \left[ 
    \begin{array}{cccc} 
     O_n &  O_n & I_n &  O_n \\ 
     O_n &  O_n & O_n & -I_n \\ 
     I_n &  O_n & O_n &  O_n \\ 
     O_n & -I_n & O_n &  O_n 
      \end{array} 
  \right] . 
\end{split} 
\label{Jposneg}
\end{equation}} 
\end{itemize}  
Let $J$ be an $\varepsilon$-paracomplex structure of $E$. 
Then such a frame field as $e$ is called 
an \textit{admissible frame field\/} of $J$. 
For an admissible frame field $e$ of $J$, 
let $\omega =[\omega^i_j ]$ be the connection form of $\nabla$ 
with respect to $e$: $\nabla e=e\omega$. 
Then we have 
\begin{equation} 
 \omega^i_j 
=\left\{ 
 \begin{array}{rl} 
 -\omega^j_i &  ((i, j) 
                \in  (\Sigma_{1   \sim 2n} \times \Sigma_{1   \sim 2n} ) 
                \cup (\Sigma_{2n+1\sim 4n} \times \Sigma_{2n+1\sim 4n} )), 
                \\ 
  \omega^j_i &  ((i, j) 
                \in  (\Sigma_{1   \sim 2n} \times \Sigma_{2n+1\sim 4n} ) 
                \cup (\Sigma_{2n+1\sim 4n} \times \Sigma_{1   \sim 2n} )), 
   \end{array} 
 \right. 
\label{omegaijgen} 
\end{equation} 
where $\Sigma_{1   \sim 2n} :=\{    1, \dots , 2n\}$, 
      $\Sigma_{2n+1\sim 4n} :=\{ 2n+1, \dots , 4n\}$. 
For $\mu \in \{ +, -\}$, we set 
$$e(\mu ):=(e_1 , \dots , e_{2n} , \mu e_{2n+1} , \dots , \mu e_{4n} ).$$ 
Suppose that $\nabla J$ is locally represented as the tensor product 
of a $1$-form $\alpha$ and an $\varepsilon$-nilpotent structure $N$ 
so that $e(\mu )$ is an admissible frame field of $N$. 
In the following, 
we say that such an $\varepsilon$-nilpotent structure as $N$ is 
\textit{related to\/} $J$ (\textit{by\/} $(e, \mu )$). 
Then using $(\nabla J)(e_k )=\nabla (Je_k )-J(\nabla e_k )$, 
we obtain 
\begin{equation} 
\begin{split} 
&   \mu                         (\omega^{n+i}_j      
  +(\varepsilon 1)^{\delta_{i1}} \omega^{3n+i}_{2n+j} )
  = \omega^{n+i}_{2n+j} 
  +(\varepsilon 1)^{\delta_{i1}} \omega^{3n+i}_j  
  = \alpha \delta^i_j , \\ 
&   \omega^i_j          -\omega^{2n+i}_{2n+j} 
  = \omega^i_{2n+j}     -\omega^{2n+i}_j 
  = 0, \\ 
&   \omega^{n+i}_{n+j}  
  -(\varepsilon 1)^{\delta_{i1} +\delta_{j1}} \omega^{3n+i}_{3n+j} 
  = \omega^{n+i}_{3n+j} 
  -(\varepsilon 1)^{\delta_{i1} +\delta_{j1}} \omega^{3n+i}_{n+j}  
  = 0     
\end{split} 
\label{omega}
\end{equation}
for $i, j\in \{1, \dots , n\}$. 
Conversely, if we suppose \eqref{omega}, 
then $\nabla J$ is locally represented as above. 
Therefore we obtain 

\begin{thm}\label{thm:hatnabla0g}
Let $J$ be an $\varepsilon$-paracomplex structure 
of $E$ for $\varepsilon \in \{ +, -\}$ 
and $e$ an admissible frame field of $J$. 
Then for $\mu \in \{ +, -\}$, the following are equivalent$:$ 
\begin{itemize} 
\item[{\rm (a)}]{$\nabla J$ is locally represented by the tensor product 
of a $1$-form $\alpha$ and an $\varepsilon$-nilpotent structure $N$ 
related to $J$ by $(e, \mu );$} 
\item[{\rm (b)}]{the connection form $\omega =[\omega^i_j ]$ of $\nabla$ 
with respect to $e$ satisfies \eqref{omega}.} 
\end{itemize} 
\end{thm}

\begin{rem}
The main object of study in the present paper is 
an $\varepsilon$-paracomplex structure $J$ of $E$ 
such that $\nabla J$ is locally represented by 
the tensor product of a $1$-form and 
an $\varepsilon$-nilpotent structure related to $J$. 
In \cite{ando8}, 
the author studies nilpotent structures of an oriented vector bundle $E$ 
of rank $4n$ with a neutral metric $h$ and an $h$-connection $\nabla$. 
For a Lie subgroup $K$ of $SO(2n, 2n)$, 
$K$-nilpotent structures of $(E, h, \nabla )$ are defined. 
For each $K$-nilpotent structure, 
a principal $K$-bundle $P$ is constructed, 
by choosing special admissible frame fields. 
In addition, $\nabla$ gives a connection in $P$, 
so that the connection form of $\nabla$ with respect to 
such an admissible frame field is valued in the Lie algebra of $K$. 
A $G$-nilpotent structure $N$ of $(E, h, \nabla )$ is characterized 
by $\nabla N=0$ (\cite{ando8}). 
If there exist a complex structure $I$ and 
paracomplex structures $J_1$, $J_2$ 
such that $h$, $\nabla$, $I$, $J_1$, $J_2$ form 
a neutral hyperK\"{a}hler structure of $E$, 
then there exist $H$-nilpotent structures of $(E, h, \nabla )$ (\cite{ando8}), 
where $H$ is as in \eqref{H}. 
In addition, 
if there exists an $H$-nilpotent structure of $(E, h, \nabla )$, 
then $N$ has the dual $H$-nilpotent structure $N^{\times}$, 
and $h$, $\nabla$, $I:=(1/2)(N+N^{\times} )$, $J_1 :=-IJ_2$, 
$J_2 :=(1/2)(N-N^{\times} )$ 
form a neutral hyperK\"{a}hler structure of $E$ (\cite{ando8}). 
\end{rem} 

In the following, suppose $n=1$. 
Let $\hat{\nabla}$ be the connection of $\bigwedge^2\!E$ induced 
by $\nabla$. 
Then $\hat{\nabla}$ is an $\hat{h}$-connection 
and it gives connections of $\bigwedge^2_+\!E$, $\bigwedge^2_-\!E$. 
Fiber bundles 
$$U_-\!\left(\textstyle\bigwedge^2_{\pm}\!E\right) 
     :=\left\{ \left. \theta \in \textstyle\bigwedge^2_{\pm}\!E 
     \ \right| \ \hat{h} (\theta , \theta )=-1\right\}$$ 
in $\bigwedge^2_{\pm}\!E$ respectively are the time-like twistor spaces 
associated with $E$. 
There exists a one-to-one correspondence between the set of sections 
of             $U_-\!\left(\textstyle\bigwedge^2_{\varepsilon}\!E\right)$ 
and the set of $\varepsilon$-paracomplex structures of $E$. 
Let $\Omega$ be a section 
of $U_-\!\left(\bigwedge^2_{\varepsilon}\!E\right)$ 
for $\varepsilon \in \{ +, -\}$. 
Then we can find $e=(e_1 , e_2 , e_3 , e_4 )$ 
satisfying $\Omega =\Omega_{\varepsilon , 2}$. 
Let $\omega =[\omega^i_j ]$ be the connection form of $\nabla$ 
with respect to $e$. 
Then we have 
\begin{equation} 
 \omega^i_j 
=\left\{ \begin{array}{rl} 
            0         &  (i=j), \\ 
          -\omega^j_i &  (\{ i, j\} \in \{ \{ 1, 2\}, \{ 3, 4\} \} ), \\ 
           \omega^j_i &  (\{ i, j\} \in \{ \{ 1, 3\}, \{ 1, 4\}, 
                                           \{ 2, 3\}, \{ 2, 4\} \} ). 
           \end{array} 
  \right. 
\label{omegaij} 
\end{equation} 
Suppose that the covariant derivative of $\Omega$ is fully light-like. 
Then for $\mu =+$ or $-$, 
we obtain 
\begin{equation} 
      \omega^3_2 +\varepsilon \omega^4_1 
=\mu (\omega^4_3 +\varepsilon \omega^2_1 ) 
\label{32+41} 
\end{equation} 
and $\hat{\nabla} \Omega =\alpha \otimes \Omega_0$, where 
\begin{equation*} 
\alpha   :=\omega^3_2 +\varepsilon \omega^4_1 , \quad 
\Omega_0 :=\Omega_{-\varepsilon , 1} +\mu \Omega_{\varepsilon , 3} . 
\label{nablaOmega+} 
\end{equation*} 
We have $\alpha \not= 0$ 
and $\Omega_0$ is a local section of 
a light-like twistor 
space $U_0\!\left(\bigwedge^2_{\varepsilon}\!E\right)$. 
Let $J$ be an $\varepsilon$-paracomplex structure of $E$ 
corresponding to $\Omega$. 
Then we obtain $Je_1 = e_3$, $Je_2 =-\varepsilon e_4$ and 
therefore $e$ is an admissible frame field of $J$. 
In addition, we obtain 
\begin{equation} 
 \nabla J 
=\alpha \otimes (\theta^1 \otimes \xi_2 -\theta^2 \otimes \xi_1 ), 
\label{nablaJ} 
\end{equation} 
where 
\begin{equation} 
\begin{array}{lcl}
\xi_1    :=e_1 -\varepsilon \mu e_3 , & \ &  
\xi_2    :=e_2 +            \mu e_4 , \\ 
\theta^1 :=e^1 +\varepsilon \mu e^3 , & \ &  
\theta^2 :=e^2 -            \mu e^4 
\end{array} 
\label{xitheta} 
\end{equation} 
and $(e^1 , e^2 , e^3 , e^4 )$ is the dual frame field 
of $(e_1 , e_2 , e_3 , e_4 )$. 
We see that $N:=\theta^1 \otimes \xi_2 -\theta^2 \otimes \xi_1$ 
is an $\varepsilon$-nilpotent structure 
corresponding to a local light-like section $\Omega_0$ and 
related to $J$ by $(e, \varepsilon \mu )$. 
Noticing that $e=(e_1 , e_2 , e_3 , e_4 )$ 
satisfying $\Omega =\Omega_{\varepsilon , 2}$ is not unique, 
we see that $\Omega_0$ and $N$ are not uniquely determined by $\Omega$. 
However, $\Omega_0$ is uniquely determined by $\Omega$ 
up to a nowhere zero function. 
Therefore $E$ has a light-like subbundle of rank 2 
given by $\Omega$ and locally generated by $\xi_1$, $\xi_2$. 
Hence we obtain 

\begin{pro}\label{pro:hatnablaO} 
Let $\Omega$ be a section 
of $U_-\!\left(\bigwedge^2_{\varepsilon}\!E\right)$ 
for $\varepsilon \in \{ +, -\}$ 
such that $\hat{\nabla} \Omega$ is fully light-like. 
Then $\hat{\nabla} \Omega$ is locally represented 
as $\hat{\nabla} \Omega =\alpha \otimes \Omega_0$ 
for a nowhere zero $1$-form $\alpha$ and a light-like section $\Omega_0$. 
In addition, if $J$ is an $\varepsilon$-paracomplex structure of $E$ 
corresponding to $\Omega$, 
then $\nabla J$ is locally represented as $\nabla J =\alpha \otimes N$, 
where $N$ is an $\varepsilon$-nilpotent structure corresponding 
to $\Omega_0$ 
and related to $J$. 
\end{pro} 

In the next section, 
we will see that the lifts of the conformal Gauss maps 
of time-like minimal surfaces in $E^3_1$ are 
examples of sections of the time-like twistor spaces 
with fully light-like covariant derivatives. 

Let $J$ be an $\varepsilon$-paracomplex structure of $E$ 
for $\varepsilon \in \{ +, -\}$. 
Suppose that the covariant derivative of a section $\Omega$ of 
$U_-\!\left(\bigwedge^2_{\varepsilon}\!E\right)$ 
corresponding to $J$ is fully light-like. 
Then whether we can find 
an ordered pseudo-orthonormal 
local frame field $e=(e_1 , e_2 , e_3 , e_4 )$ of $E$ 
as above for $\Omega$ 
such that $\Omega_0$ is horizontal is determined by $J$. 
Referring to \cite{ando7}, 
we see that $\Omega_0$ is horizontal if and only if 
the corresponding nilpotent structure $N$ satisfies $\nabla N=0$. 
The following proposition will be used in Section~\ref{sect:flatdim4}. 

\begin{pro}\label{pro:O-1O+3hor} 
The following two conditions are equivalent to each other\/$:$ 
\begin{itemize} 
\item[{\rm (a)}]{there exists 
an ordered pseudo-orthonormal local frame field $e$ 
of $E$ as above for $\Omega$ 
such that $\Omega_0$ is horizontal\/$;$} 
\item[{\rm (b)}]{there exists 
an ordered pseudo-orthonormal local frame field $e$ 
of $E$ as above for $\Omega$ 
such that the connection form $[\omega^i_j ]$ of $\nabla$ 
with respect to $e$ 
satisfies $d(\omega^4_2 -\varepsilon \omega^3_1 )=0$.} 
\end{itemize} 
\end{pro} 

\vspace{3mm} 

\par\noindent 
\textit{Proof} \ 
Suppose (a). 
Then we see that $\Omega_0$ is horizontal 
for an ordered pseudo-orthonormal 
local frame field $e$. 
Using \eqref{32+41}, we obtain 
\begin{equation} 
 \sqrt{2} \hat{\nabla} \Omega_0 
=\hat{\nabla} (\xi_1 \wedge \xi_2 ) 
=\mu (\omega^4_2 -\varepsilon \omega^3_1 )\otimes \xi_1 \wedge \xi_2 . 
\label{42-+31} 
\end{equation} 
Therefore $\hat{\nabla} \Omega_0 =0$ 
means $\omega^4_2 -\varepsilon \omega^3_1 =0$. 
Hence $e$ is a suitable frame field for Condition (b). 
Suppose (b) and let $e$ be as in (b). 
Then there exists a function $f$ 
satisfying $\omega^4_2 -\varepsilon \omega^3_1 =-df$. 
Using this and \eqref{42-+31}, 
we obtain $\hat{\nabla} (e^{\mu f} \xi_1 \wedge \xi_2 )=0$. 
We set 
\begin{equation*} 
\begin{array}{lcl} 
 e'_1 :=e_1 , & \ & 
 e'_2 := ({\rm cosh}\,f)e_2 +({\rm sinh}\,f)e_4 , \\ 
 e'_3 :=e_3 , & \ & 
 e'_4 := ({\rm sinh}\,f)e_2 +({\rm cosh}\,f)e_4 . 
  \end{array} 
\end{equation*} 
Then $e'=(e'_1 , e'_2 , e'_3 , e'_4 )$ is 
an ordered pseudo-orthonormal local frame field satisfying 
\begin{itemize} 
\item[{\rm (i)}]{$e'_1$, $e'_2$ are space-like 
and   $e'_3$, $e'_4$ are time-like,} 
\item[{\rm (ii)}]{$e'$ gives the orientation of $E$,} 
\item[{\rm (iii)}]{$\Omega =\Omega'_{\varepsilon , 2}$, 
where $\Omega'_{\pm , 2}$ are defined for $e'$ 
as in \eqref{Omega+-123},} 
\item[{\rm (iv)}]{$\Omega_0$ for $e'$ is horizontal.} 
\end{itemize} 
Hence $e'=(e'_1 , e'_2 , e'_3 , e'_4 )$ is a suitable frame field 
for Condition (a). 
\hfill 
$\square$ 

\section{The conformal Gauss maps of time-like minimal surfaces 
in \mbox{\boldmath{$E^3_1$}}}\label{sect:cgmtms} 

\setcounter{equation}{0} 

Let $M$ be a Lorentz surface. 
Let $\iota : M\longrightarrow E^3_1$ be a time-like and conformal 
immersion of $M$ into $E^3_1$. 
Let $\langle \ , \ \rangle$ be the metric of $E^5_2$ and set 
$$L:=\{ x=(x^1 , x^2 , x^3 , x^4 , x^5 )\in E^5_2 \setminus \{ 0\} \ | \ 
        \langle x, x\rangle =0\} .$$ 
We identify $E^3_1$ with $L\cap \{ x^5 =x^1 +1\}$ and therefore 
we consider $\iota$ to be an $L$-valued function. 
Suppose that $\iota$ is minimal. 
Let $\gamma$ be the unit normal vector field of $\iota$ 
determined by the orientations of $M$ and $E^3_1$. 
Then $\gamma$ is a map from $M$ into $S^4_2$. 
Let $g^M$ be the induced metric of $M$ by $\iota$ 
and $K^M$ the curvature of $g^M$. 
Let ${\rm Reg}\,(\iota )$ be the set of nonzero points of $K^M$. 
Then the restriction of $\gamma$ on ${\rm Reg}\,(\iota )$ is 
a time-like and conformal immersion which induces a Lorentz metric $g$ 
given by $g=-K^M g^M$. 
We call $\gamma : M\longrightarrow S^4_2$ 
the \textit{conformal Gauss map\/} of $\iota$. 
We can refer to \cite{AHHK} 
for the conformal Gauss maps of time-like minimal surfaces in $E^3_1$. 

Let $w=u+jv$ be a local paracomplex coordinate of ${\rm Reg}\,(\iota )$. 
Then we represent $g^M$ as $g^M =e^{2\lambda} dwd\overline{w}$ 
and we see that there exist functions $l$, $m$ satisfying 
\begin{equation} 
\iota_{uu} =\lambda_u \iota_u +\lambda_v \iota_v +l\gamma , \quad 
\iota_{uv} =\lambda_v \iota_u +\lambda_u \iota_v +m\gamma 
\label{iuuiuv} 
\end{equation} 
and $\iota_{vv} =\iota_{uu}$. 
By \eqref{iuuiuv}, we obtain 
\begin{equation} 
 (\gamma_u \ \gamma_v ) 
=(\iota_u  \ \iota_v  ) 
  \dfrac{1}{e^{2\lambda}} 
  \left[ \begin{array}{cc} 
         -l & -m \\ 
          m &  l 
          \end{array} 
 \right] . 
\label{gugv} 
\end{equation} 
Therefore we can consider $\iota$ to be a light-like normal vector field 
of $F:=\gamma |_{{\rm Reg}\,(\iota )}$. 
Since $K^M \not= 0$ on ${\rm Reg}\,(\iota )$, we have $l^2 \not= m^2$. 
Let $A_{\iota}$ be the shape operator of $F$ with respect to $\iota$. 
Then from \eqref{gugv}, we obtain 
\begin{equation} 
\begin{split} 
&  \left( A_{\iota} \left( \dfrac{\partial}{\partial u} \right) \ 
          A_{\iota} \left( \dfrac{\partial}{\partial v} \right) 
   \right) \\ 
& =\left( dF\left( \dfrac{\partial}{\partial u} \right) \ 
          dF\left( \dfrac{\partial}{\partial v} \right) 
   \right) 
   \dfrac{e^{2\lambda}}{l^2 -m^2} 
   \left[ \begin{array}{cc} 
           l &  m \\ 
          -m & -l 
            \end{array} 
   \right] . 
\end{split} 
\label{Aiota} 
\end{equation} 
Let $\nu$ be a light-like normal vector field of $F$ 
satisfying $\langle \nu , \iota \rangle =-1$. 

\begin{pro}\label{pro:nu} 
The shape operator $A_{\nu}$ of $F$ with respect to $\nu$ vanishes. 
\end{pro} 

\vspace{3mm} 

\par\noindent 
\textit{Proof} \ 
We can prove Proposition~\ref{pro:nu}, referring to \cite{ando4}. 
Since $\iota$ is minimal, 
a paracomplex quartic differential $\Tilde{Q}$ defined on $M$ vanishes. 
This means that a paracomplex quartic differential $Q$ 
defined on ${\rm Reg}\,(\iota )$ and given by $F$ vanishes. 
Therefore $\nu$ is contained in a constant direction in $E^5_2$. 
This means that $A_{\nu}$ vanishes. 
\hfill 
$\square$ 

\vspace{3mm} 

From \eqref{Aiota} and Proposition~\ref{pro:nu}, 
we see that $F$ has zero mean curvature vector. 
This is also seen from a fact that 
a time-like minimal immersion $\iota$ is of Willmore type (\cite{ando4}). 

Suppose $K^M <0$, i.e., $l^2 >m^2$. 
Then $F_u$ is space-like and $F_v$ is time-like. 
We set 
$$e_1 :=\dfrac{1}{\sqrt{-K^M} e^{\lambda}} F_u , \quad 
  e_3 :=\dfrac{1}{\sqrt{-K^M} e^{\lambda}} F_v .$$ 
Let $e_2$, $e_4$ be normal vector fields of $F$ satisfying 
$$\langle e_2 , e_2 \rangle = 1, \ 
  \langle e_4 , e_4 \rangle =-1, \ 
  \iota =\dfrac{1}{\sqrt{2}} (e_4 -e_2 ), \ 
  \nu   =\dfrac{1}{\sqrt{2}} (e_4 +e_2 ).$$ 
Then $e=(e_1 , e_2 , e_3 , e_4 )$ is 
an ordered pseudo-orthonormal local frame field 
of the pull-back bundle $E:= F^*\!T\!S^4_2$ 
by $F:{\rm Reg}\,(\iota )\longrightarrow S^4_2$ 
and we can suppose that $e=(e_1 , e_2 , e_3 , e_4 )$ gives the orientation 
of $S^4_2$. 
Let $(e^1 , e^3 )$ be 
the dual frame field of an ordered pseudo-orthonormal local 
frame field $(e_1 , e_3 )$ of the tangent bundle $T\!M$. 
Let $h$, $\nabla$ be the metric and the Levi-Civita connection of $S^4_2$. 
Then they naturally give a metric and an $h$-connection of $E$ 
and these are also denoted by $h$, $\nabla$. 
Let $\omega =[\omega^i_j ]$ be the connection form of $\nabla$ 
with respect to $e$. 
Then we obtain 
\begin{equation} 
\omega^2_1 =\omega^4_1 =-(\Tilde{l} e^1 +\Tilde{m} e^3 ), \quad 
\omega^3_2 =\omega^4_3 =-(\Tilde{m} e^1 +\Tilde{l} e^3 )
\label{omegacgm} 
\end{equation} 
with 
$$\Tilde{l} :=\dfrac{le^{2\lambda}}{\sqrt{2} (l^2 -m^2 )} , \quad 
  \Tilde{m} :=\dfrac{me^{2\lambda}}{\sqrt{2} (l^2 -m^2 )} .$$ 
Let $\hat{\nabla}$ be the connection of $\bigwedge^2 E$ 
induced by $\nabla$. 
Then by \eqref{omegacgm}, 
we obtain $\hat{\nabla} \Omega_{\varepsilon , 2} =\alpha \otimes \Omega_0$ 
for $\varepsilon \in \{ +, -\}$, 
where 
$$\alpha   :=-(\varepsilon \Tilde{l} +\Tilde{m} )(e^1 +\varepsilon e^3 ), 
  \quad 
  \Omega_0 :=\Omega_{-\varepsilon , 1} +\Omega_{\varepsilon , 3}$$ 
and $\Omega_{\pm , 1}$, $\Omega_{\pm , 2}$, $\Omega_{\pm , 3}$ are 
as in \eqref{Omega+-123}. 
We have $\alpha \not= 0$. 
Let $\hat{F}_{\pm}$ be the lifts of $F$. 
Then $\hat{F}_{\pm}$ are sections 
of $U_-\!\left(\bigwedge^2_{\pm}\!E\right)$ 
and locally represented as $\hat{F}_{\pm} =\Omega_{\pm , 2}$ respectively. 
Let $J_{\varepsilon}$ be an $\varepsilon$-paracomplex structure of $E$ 
corresponding to $\hat{F}_{\varepsilon}$. 
Then we have $J_{\varepsilon} e_1 =             e_3$, 
             $J_{\varepsilon} e_2 =-\varepsilon e_4$. 
For $J=J_{\varepsilon}$, we obtain \eqref{nablaJ} 
with \eqref{xitheta} and $\mu =+$. 
We see from \eqref{nablaJ} for $J=J_{\varepsilon}$ 
that $\nabla J_{\varepsilon}$ is valued in 
a light-like subbundle of $E=F^*\!T\!S^4_2$ 
of rank $2$ 
which is locally generated by $e_1 -\varepsilon e_3$, $e_2 +e_4$. 
Hence we obtain 

\begin{thm}\label{thm:cgmtms} 
The covariant derivatives of the lifts $\hat{F}_{\pm}$ are 
fully light-like. 
In addition, for an $\varepsilon$-paracomplex structure $J_{\varepsilon}$ 
corresponding to $\hat{F}_{\varepsilon}$, 
$\nabla J_{\varepsilon}$ is locally represented 
as $\nabla J_{\varepsilon} =\alpha \otimes N$, 
where $\alpha :=-(\varepsilon \Tilde{l} +\Tilde{m} )\theta^1$, 
      $N      := \theta^1 \otimes \xi_2 -\theta^2 \otimes \xi_1$ 
with \eqref{xitheta} and $\mu =+$. 
\end{thm} 

Let $M$ be a Lorentz surface 
and $F:M\longrightarrow S^4_2$ a time-like and conformal immersion 
with zero mean curvature vector. 
According to \cite{ando5}, 
if the covariant derivatives of the lifts $\hat{F}_{\pm}$ are 
fully light-like, 
then $F$ has one of the following two properties: 
\begin{itemize} 
\item[{\rm (i)}]{the shape operator of $F$ 
with respect to a light-like normal vector field 
vanishes\/$;$} 
\item[{\rm (ii)}]{the shape operator of $F$ 
with respect to any normal vector field is 
light-like or zero.} 
\end{itemize} 
As was already seen by Proposition~\ref{pro:nu}, 
the conformal Gauss maps of time-like minimal surfaces in $E^3_1$ 
satisfy (i). 
We can find a characterization of 
a time-like and conformal immersion $F:M\longrightarrow S^4_2$ 
with zero mean curvature vector and (ii) in terms of 
the Gauss-Codazzi-Ricci equations (\cite{ando5}). 
The immersion $F$ satisfies either (i) or (ii) 
if and only if the curvature $K$ of the induced metric on $M$ by $F$ is 
identically equal to $1$ (\cite{ando6}). 
In addition, $K\equiv 1$ if and only if 
not only the paracomplex quartic differential $Q$ is null or zero 
but also the normal connection $\nabla^{\perp}$ of $F$ is flat (\cite{ando6}). 
Notice that $Q$ is null or zero if and only if 
at least one of the covariant derivatives of $\hat{F}_{\pm}$ is 
light-like or zero (\cite{ando6}) and 
that $K\equiv 1$ means that $\nabla^{\perp}$ is flat, 
while the converse is not necessarily true (\cite{ando6}). 

\section{Walker \mbox{\boldmath{$4n$}}-manifolds}\label{sect:wm} 

\setcounter{equation}{0} 

Let $M$ be an oriented neutral $4n$-manifold and $h$ its metric. 
Let $\nabla$ be the Levi-Civita connection of $h$. 
Let $\mathcal{D}$ be a light-like $2n$-dimensional distribution on $M$. 
We say that $(M, h, \mathcal{D} )$ is a \textit{Walker manifold\/} 
if the covariant derivatives of local generators of $\mathcal{D}$ 
by any tangent vector of $M$ with respect to $\nabla$ are contained 
in $\mathcal{D}$. 
See \cite{BGGNV}, \cite{walker} for Walker manifolds. 
If $(M, h, \mathcal{D} )$ is a Walker manifold, 
then $\mathcal{D}$ is involutive. 
In the case of $n=1$, $\mathcal{D}$ is involutive 
if and only if the covariant derivatives of local generators 
of $\mathcal{D}$ by any tangent vector in $\mathcal{D}$ are 
contained in $\mathcal{D}$ (see \cite{ando7}). 

Let $J$ be an almost $\varepsilon$-paracomplex structure of $M$ 
for $\varepsilon \in \{ +, -\}$. 
Suppose that $\nabla J$ satisfies (a) in Theorem~\ref{thm:hatnabla0g} 
for a nowhere zero $1$-form $\alpha$. 
Then $\nabla J$ gives a light-like subbundle of $T\!M$ of rank $2n$, 
i.e., a light-like $2n$-dimensional distribution $\mathcal{D}_J$ on $M$. 

\begin{thm}\label{thm:JDW} 
Let $J$ be as above. 
Then $(M, h, \mathcal{D}_J )$ is a Walker manifold. 
\end{thm} 

\vspace{3mm} 

\par\noindent 
\textit{Proof} \ 
Let $e=(e_1 ,\dots e_{2n} , e_{2n+1} , \dots , e_{4n} )$ be 
an admissible frame field of $J$ 
such that $\nabla J$ is locally represented as $\nabla J=\alpha \otimes N$ 
for an almost $\varepsilon$-nilpotent structure $N$ 
related to $J$ by $(e, \mu )$. 
Then we have \eqref{omega}. 
Let $\xi_1 , \dots , \xi_{2n}$ be as in \eqref{xi2} for $e(\mu )$. 
For $i, j\in \{ 1, \dots , n\}$, we obtain 
\begin{equation} 
\begin{split} 
&   h(\nabla \xi_i , \xi_j ) 
  = \omega^j_i -\mu \omega^j_{2n+i} 
               +\mu \omega^{2n+j}_i -\omega^{2n+j}_{2n+i} , \\ 
&   h(\nabla \xi_i , \xi_{n+j} ) 
  =-h(\xi_i , \nabla \xi_{n+j} ) \\ 
& = \omega^{n+j}_i  
   -\mu                               \omega^{n+j}_{2n+i} 
   -\mu (\varepsilon 1)^{\delta_{j1}} \omega^{3n+j}_i 
   +    (\varepsilon 1)^{\delta_{j1}} \omega^{3n+j}_{2n+i} , \\ 
&   h(\nabla \xi_{n+i} , \xi_{n+j} ) \\ 
& =                                                \omega^{n+j}_{n+i}  
   +\mu (\varepsilon 1)^{\delta_{i1}}              \omega^{n+j}_{3n+i} 
   -\mu (\varepsilon 1)^{\delta_{j1}}              \omega^{3n+j}_{n+i} 
   -    (\varepsilon 1)^{\delta_{i1} +\delta_{j1}} \omega^{3n+j}_{3n+i} . 
\end{split} 
\label{nablaxiixij} 
\end{equation} 
Applying \eqref{omega} to \eqref{nablaxiixij}, 
we obtain $h(\nabla \xi_i , \xi_j )=0$ for $i, j=1, \dots , 2n$. 
Therefore $(M, h, \mathcal{D}_J )$ is a Walker manifold. 
\hfill 
$\square$ 

\vspace{3mm} 

\begin{rem}\label{rmk:rmk2}  
Suppose $n=1$. 
Let $J$ be an almost $\varepsilon$-paracomplex structure of $M$ 
such that the covariant derivative of the corresponding section 
of $U_-\!\left(\textstyle\bigwedge^2_{\varepsilon}\!T\!M\right)$ is 
fully light-like. 
Then noticing Proposition~\ref{pro:hatnablaO}, 
we see from Theorem~\ref{thm:JDW} 
that $(M, h, \mathcal{D}_J )$ is a Walker manifold. 
\end{rem} 

Suppose $n=1$. 
Let $\mathcal{D}$ be a light-like two-dimensional distribution on $M$. 
Then $\mathcal{D}$ is locally generated by 
light-like vector fields $\xi_1$, $\xi_2$ 
such that $(1/\sqrt{2})\xi_1 \wedge \xi_2$ is a local section 
of $U_0\!\left(\textstyle\bigwedge^2_{\varepsilon}\!T\!M\right)$ 
for $\varepsilon =+$ or $-$. 
Therefore there exists an almost $\varepsilon$-nilpotent structure $N$ 
of a neighborhood of each point of $M$ 
such that $\mathcal{D}$ is locally given by $\pi_N$. 
We will prove 

\begin{thm}\label{thm:JDW2} 
Suppose $n=1$. 
Let $\mathcal{D}$ be a light-like two-dimensional distribution on $M$ 
such that $(M, h, \mathcal{D} )$ is a Walker manifold. 
Let $N$ be an almost $\varepsilon$-nilpotent structure as above. 
Then there exists an almost $\varepsilon$-paracomplex structure $J$ 
of a neighborhood of the point satisfying 
\begin{itemize} 
\item[{\rm (a)}]{$\nabla J=\alpha \otimes N$ 
for a nowhere zero $1$-form $\alpha$,} 
\item[{\rm (b)}]{$N$ is related to $J$.} 
\end{itemize} 
\end{thm} 

\vspace{3mm} 

\par\noindent 
\textit{Proof} \ 
Let $\mathcal{D}$, $N$ be as in Theorem~\ref{thm:JDW2}. 
Suppose that $N$ corresponds to a local section $\Omega_0$ 
of $U_0\!\left(\bigwedge^2_+\!T\!M\right)$. 
Let $(e_1 , e_2 , e_3 , e_4 )$ be an admissible frame field of $N$. 
Then $\Omega_0 =\Omega_{-, 1} +\Omega_{+, 3}$ 
and $\mathcal{D}$ is locally generated 
by $\xi_1 =e_1 -e_3$, $\xi_2 =e_2 +e_4$. 
Since $(M, h, \mathcal{D} )$ is a Walker manifold, 
we have $h(\nabla \xi_1 , \xi_2 )=0$ and 
this is rewritten into $\omega^3_2 +\omega^4_1 =\omega^4_3 +\omega^2_1$. 
In addition, we can suppose $\omega^3_2 +\omega^4_1 \not= 0$, 
rechoosing $(e_1 , e_2 , e_3 , e_4 )$ if necessary: 
if we set 
$$(\Tilde{e}_1 , \Tilde{e}_2 , \Tilde{e}_3 , \Tilde{e}_4 ) 
 =(e_1  , e_2  , e_3  , e_4 )
  \left[ \begin{array}{cccc} 
         1 &  0               & 0 &   0               \\ 
         0 & \frac{h^2 +2}{2} & h & -\frac{h^2}{2}    \\ 
         0 &       h          & 1 & -      h          \\ 
         0 & \frac{h^2}{2}    & h & -\frac{h^2 -2}{2} 
           \end{array} 
  \right]$$  
for a function $h$, 
then 
\begin{itemize} 
\item[{\rm (i)}]{$(\Tilde{e}_1 , \Tilde{e}_2 , \Tilde{e}_3 , \Tilde{e}_4 )$ is 
an ordered pseudo-orthonormal local frame field of $M$ 
giving the orientation of $M$,} 
\item[{\rm (ii)}]{$(\Tilde{e}_1 , \Tilde{e}_2 , \Tilde{e}_3 , \Tilde{e}_4 )$ 
is an admissible frame field of $N$,} 
\item[{\rm (iii)}]{for the connection form $[\Tilde{\omega}^i_j ]$ 
of $\nabla$ with respect to 
$(\Tilde{e}_1 , \Tilde{e}_2 , \Tilde{e}_3 , \Tilde{e}_4 )$, 
$$\Tilde{\omega}^3_2 +\Tilde{\omega}^4_1 
 =dh -\dfrac{h^2 -2}{2} (\omega^3_2 +\omega^4_1 ) 
     +\dfrac{h^2}{2}    (\omega^4_3 +\omega^2_1 )  
     +       h          (\omega^4_2 -\omega^3_1 ).$$} 
\end{itemize} 
Therefore the covariant derivative of $\Omega =\Omega_{+, 2}$ is 
fully light-like. 
Let $J$ be an almost $+$-paracomplex structure 
corresponding to $\Omega$. 
Then $J$ satisfies (a), (b) in Theorem~\ref{thm:JDW2}. 
In the case where $\Omega_0$ is a local section 
of $U_0\!\left(\bigwedge^2_-\!T\!M\right)$, 
we can obtain the same result. 
\hfill 
$\square$ 

\vspace{3mm} 

\begin{rem} 
Let $N$ be an almost $\varepsilon$-nilpotent structure of $M$ 
satisfying $\nabla N=0$. 
Then $(M, h, \mathcal{D} )$ is a Walker manifold 
for a light-like two-dimensional distribution $\mathcal{D}$ on $M$ 
given by $\pi_N$. 
Therefore there exists 
an almost $\varepsilon$-paracomplex structure $J$ 
on a neighborhood of each point of $M$ satisfying (a), (b) 
in Theorem~\ref{thm:JDW2}. 
\end{rem} 

\section{The square norm of the covariant derivative of 
an \mbox{\boldmath{$h$}}-reversing almost paracomplex structure} 

\setcounter{equation}{0} 

Let $M$, $h$, $\nabla$ be as in the beginning of the previous section. 
Let $J$ be an $h$-reversing almost paracomplex structure of $M$. 
Then the square norm $\parallel\!\nabla J\!\parallel^2$ of $\nabla J$ is 
defined by 
$$\parallel\!\nabla J\!\parallel^2 
 :=\sum^{4n}_{i, j=1} \varepsilon_i \varepsilon_j 
h((\nabla_{e_i} J)(e_j ), (\nabla_{e_i} J)(e_j )),$$ 
where $e=(e_1 , \dots , e_{2n} , e_{2n+1} , \dots , e_{4n} )$ is 
an ordered pseudo-orthonormal local frame field of $M$ 
giving the orientation and 
$\varepsilon_k =1$, $\varepsilon_{2n+k} =-1$ ($k=1, \dots , 2n$). 
The square norm $\parallel\!\nabla J\!\parallel^2$ does not depend 
on the choice of $e$. 
Let $\Omega^*$ be a $2$-form on $M$ defined 
by $\Omega^* (X, Y):=h(X, JY)$ for tangent vectors $X$, $Y$ of $M$. 
Then the square norm $\parallel\!\nabla \Omega^*\!\parallel^2$ 
of $\nabla \Omega^*$ is defined by 
\begin{equation} 
  \parallel\!\nabla \Omega^*\!\parallel^2 
:=\sum^{4n}_{i, j, k=1} \varepsilon_i \varepsilon_j \varepsilon_k 
((\nabla_{e_i} \Omega^* )(e_j , e_k ))^2 . 
\label{omega*} 
\end{equation} 
Noticing 
$$h((\nabla_{e_i} J)(e_j ), e_k )=-(\nabla_{e_i} \Omega^* )(e_j , e_k ),$$ 
we obtain $\parallel\!\nabla J\!\parallel^2 
          =\parallel\!\nabla \Omega^*\!\parallel^2$. 

Suppose that $J$ is an almost $\varepsilon$-paracomplex structure of $M$. 
Let $e$ be an admissible frame field of $J$. 
Then we have 
$$\Omega^* 
 =\sum^n_{i=1} (                              e^i     \wedge e^{2n+i} 
               -(\varepsilon 1)^{\delta_{i1}} e^{n+i} \wedge e^{3n+i} ),$$ 
where $(e^1 , \dots , e^{2n} , e^{2n+1} , \dots , e^{4n} )$ is 
the dual frame field of $e$. 
Suppose that $\nabla J$ satisfies (a) in Theorem~\ref{thm:hatnabla0g}. 
Then we have \eqref{omega}. 
Therefore we see 
\begin{itemize} 
\item[{\rm (i)}]{$((\nabla_{e_i} \Omega^* )(e_j , e_k ))^2 =\alpha (e_i )^2$ 
if 
\begin{equation} 
\{ j, k\} \in \{ \{   l, n+l\} ,  \{    l, 3n+l\} , 
                 \{ n+l, 2n+l\} , \{ 2n+l, 3n+l\} \}  
\label{jk} 
\end{equation} 
for $l\in \{ 1, \dots , n\}$,} 
\item[{\rm (ii)}]{$((\nabla_{e_i} \Omega^* )(e_j , e_k ))^2 =0$ 
if $\{ j, k\}$ does not satisfy \eqref{jk} 
for any $l\in \{ 1, \dots , n\}$.} 
\end{itemize}  
Therefore noticing 
$$\varepsilon_j \varepsilon_k 
 =\left\{ \begin{array}{rl} 
           1 &  (\{ j, k\} \in \bigcup^n_{l=1} 
                           \{ \{   l,  n+l\} , \{ 2n+l, 3n+l\} \} ), \\ 
          -1 &  (\{ j, k\} \in \bigcup^n_{l=1} 
                           \{ \{   l, 3n+l\} , \{  n+l, 2n+l\} \} ), 
            \end{array} 
  \right.$$ 
we see by \eqref{omega*} 
that $\parallel\!\nabla \Omega^*\!\parallel^2 =0$. 
This means $\parallel\!\nabla J\!\parallel^2 =0$. 

Hence we obtain 

\begin{thm}\label{thm:snorm0} 
Let $J$ be an almost $\varepsilon$-paracomplex structure of $M$ 
such that $\nabla J$ is locally represented 
as in {\rm (a)} in Theorem~\ref{thm:hatnabla0g}. 
Then $\parallel\!\!\nabla J\!\!\parallel^2 =0$, that is, 
$(M, h)$ equipped with $J$ is isotropic paraK\"{a}hler. 
\end{thm} 

\section{Frame fields of neutral vector bundles 
of rank \mbox{\boldmath{$4n$}} 
with flat metric connections}\label{sect:flatdim4n}  

\setcounter{equation}{0} 

Let $E$ be an oriented vector bundle 
over $\mbox{\boldmath{$R$}}^m$ of rank $4n$ ($m\geq 2$, $n\geq 1$). 
Let $h$ be a neutral metric of $E$ and $\nabla$ an $h$-connection of $E$. 
Suppose that $\nabla$ is flat, that is, 
the curvature tensor $R$ of $\nabla$ vanishes. 

Let $J$ be a $+$-paracomplex structure of $E$. 
Then there exists an ordered pseudo-orthonormal frame field 
$e=(e_1 , \dots , e_{2n} , e_{2n+1} , \dots , e_{4n} )$ of $E$ 
giving the orientation and satisfying \eqref{Jposneg} 
with $\varepsilon =+$. 
Let $\omega =[\omega^i_j ]$ be the connection form of $\nabla$ 
with respect to $e$. 
Then we have \eqref{omegaijgen}. 
Suppose that $\nabla J$ is nowhere zero and 
satisfies (a) in Theorem~\ref{thm:hatnabla0g} for $\varepsilon =+$. 
Then we have \eqref{omega} for a nowhere zero 1-form $\alpha$ 
on $\mbox{\boldmath{$R$}}^m$. 
Suppose $\mu =+$. 
Let $D_{ij}$ be the $(i, j)$-block 
of $\omega$ ($i, j\in \{ 1, 2, 3, 4\}$). 
Then we have 
$$\omega 
 =\left[ \begin{array}{cccc} 
          D_{11} & D_{12} & D_{13} & D_{14} \\ 
          D_{21} & D_{22} & D_{23} & D_{24} \\ 
          D_{31} & D_{32} & D_{33} & D_{34} \\ 
          D_{41} & D_{42} & D_{43} & D_{44} 
           \end{array} 
  \right] .$$ 
From \eqref{omegaijgen}, we obtain 
\begin{equation} 
{}^t\!D_{ij} 
=\left\{ 
 \begin{array}{rl} 
  -D_{ji} &  (i, j=1, 2 \ {\rm or} \ i, j=3, 4), \\ 
   D_{ji} &  (i=1, 2, j=3, 4 \ {\rm or} \ i=3, 4, j=1, 2). 
   \end{array} 
 \right. 
\label{Dij} 
\end{equation} 
From \eqref{omega}, we obtain 
\begin{equation} 
\begin{split} 
& D_{11} =D_{33} , \quad 
  D_{13} =D_{31} , \quad 
  D_{22} =D_{44} , \quad 
  D_{24} =D_{42} , \\ 
& D_{21} +D_{43} =D_{41} +D_{23} =\alpha I_n . 
\end{split} 
\label{D} 
\end{equation} 
Let $\Psi_{ij}$ be the $(i, j)$-block 
of $\omega \wedge \omega$ ($i, j\in \{ 1, 2, 3, 4\}$). 
Suppose 
\begin{equation} 
D_{11} = D_{22} , \quad 
D_{13} = D_{24} . 
\label{D2} 
\end{equation} 
Then from \eqref{D} and \eqref{D2}, 
we obtain $\Psi_{21} +\Psi_{43} =O_n$, $\Psi_{41} +\Psi_{23} =O_n$. 
Since $R$ vanishes, $d\omega +\omega \wedge \omega$ vanishes. 
Therefore we obtain $d(D_{21} +D_{43} )=O_n$, $d(D_{41} +D_{23} )=O_n$. 
These mean $d\alpha =0$. 
Therefore a function $f$ on $\mbox{\boldmath{$R$}}^m$ 
satisfies $df=\alpha$. 
Suppose $D_{21} =D_{43}$, $D_{41} =D_{23}$. 
Then we see by \eqref{Dij} that $\omega$ is represented as in the form of 
\begin{equation} 
  \omega 
 =\left[ \begin{array}{cccc} 
           D_{11}            & -\frac{1}{2} df I_n & 
           D_{31}            &  \frac{1}{2} df I_n \\ 
          \frac{1}{2} df I_n &   D_{11}            & 
          \frac{1}{2} df I_n &   D_{31}            \\ 
           D_{31}            &  \frac{1}{2} df I_n & 
           D_{11}            & -\frac{1}{2} df I_n \\ 
          \frac{1}{2} df I_n &   D_{31}            & 
          \frac{1}{2} df I_n &   D_{11} 
           \end{array} 
  \right] , 
\label{omegaD} 
\end{equation} 
where $D_{11}$, $D_{31}$ satisfy 
${}^t\!D_{11} =-D_{11}$, ${}^t\!D_{31} =D_{31}$ and 
\begin{equation} 
\begin{split} 
& dD_{11} +D_{11} \wedge D_{11} +D_{31} \wedge D_{31} =O_n , \\ 
& dD_{31} +D_{31} \wedge D_{11} +D_{11} \wedge D_{31} =O_n . 
\end{split} 
\label{dD} 
\end{equation} 

Let $D_{11}$, $D_{31}$ be $n\times n$ matrices 
such that each component is a $1$-form on $\mbox{\boldmath{$R$}}^m$. 
Suppose that $D_{11}$, $D_{31}$ satisfy 
${}^t\!D_{11} =-D_{11}$, ${}^t\!D_{31} =D_{31}$ and \eqref{dD}. 
Let $f$ be a function on $\mbox{\boldmath{$R$}}^m$ 
such that $df$ is nowhere zero. 
Then a $4n\times 4n$ matrix $\omega$ defined by \eqref{omegaD} 
satisfies \eqref{omegaijgen} and \eqref{omega} 
with $\alpha =df$ and $\varepsilon =\mu =+$, 
and $d\omega +\omega \wedge \omega$ vanishes. 
Therefore, noticing $R=0$, 
we see that there exists 
an ordered frame 
field $e=(e_1 , \dots , e_{2n} , e_{2n+1} , \dots , e_{4n} )$ 
of $E$ such that $\omega$ is a connection form of $\nabla$ 
with respect to $e$. 
Such a frame field is uniquely determined by an initial value at a point. 
If $e$ satisfies 
\begin{equation} 
 h(e_i , e_j ) 
=\left\{ \begin{array}{rl} 
          1 &  (i=j=1,    \dots , 2n), \\ 
         -1 &  (i=j=2n+1, \dots , 4n), \\ 
          0 &  (i\not= j) 
           \end{array} 
 \right.
\label{xipxiq4n} 
\end{equation} 
at a point of $\mbox{\boldmath{$R$}}^m$, 
then noticing that $\nabla$ is an $h$-connection, 
we see that $e$ satisfies \eqref{xipxiq4n} on $\mbox{\boldmath{$R$}}^m$. 
Therefore there exists 
an ordered pseudo-orthonormal frame field $e$ 
of $E$ giving the orientation of $E$ and satisfying $\nabla e=e\omega$. 
Let $J$ be a $+$-paracomplex structure of $E$ 
satisfying \eqref{Jposneg} for $\varepsilon =+$ and $e$. 
Then $\nabla J$ is nowhere zero and 
satisfies (a) in Theorem~\ref{thm:hatnabla0g} 
with $\alpha =df$ for $\varepsilon =\mu =+$. 

\begin{ex}\label{ex:ex1} 
We set $D_{11} :=O_n$ and $D_{31} :=d\phi C_0$, 
where $\phi$ is a function on $\mbox{\boldmath{$R$}}^m$ 
and $C_0$ is a constant symmetric matrix. 
Then $D_{11}$, $D_{31}$ satisfy 
${}^t\!D_{11} =-D_{11}$, ${}^t\!D_{31} =D_{31}$ and \eqref{dD}. 
\end{ex} 

\begin{ex}\label{ex:ex2}  
We set $D_{11} :=O_n$ and  
$$D_{31} 
:=\left[ \begin{array}{ccccc} 
           df_1  &  df_2  &   0    & \cdots &  0     \\ 
           df_2  &  df_1  &  df_2  & \ddots & \vdots \\ 
            0    &  df_2  & \ddots & \ddots &  0     \\ 
          \vdots & \ddots & \ddots & \ddots &  df_2  \\ 
            0    & \cdots &   0    &  df_2  &  df_1 
           \end{array} 
  \right]$$ 
for functions $f_1$, $f_2$ on $\mbox{\boldmath{$R$}}^m$. 
Then $D_{11}$, $D_{31}$ satisfy 
${}^t\!D_{11} =-D_{11}$, 
${}^t\!D_{31} = D_{31}$ and \eqref{dD}. 
\end{ex} 

\begin{ex}\label{ex:ex3}  
We set $D_{31} :=O_n$ and $D_{11} :=d\psi C_0$, 
where $\psi$ is a function on $\mbox{\boldmath{$R$}}^m$ 
and $C_0$ is a constant skew-symmetric matrix. 
Then $D_{11}$, $D_{31}$ satisfy 
${}^t\!D_{11} =-D_{11}$, ${}^t\!D_{31} =D_{31}$ and \eqref{dD}. 
\end{ex} 

\begin{ex}\label{ex:ex4}  
Suppose $n=4p$ for $p\in \mbox{\boldmath{$N$}}$. 
We set  
$$D_{i1} 
:=\left[ \begin{array}{ccccc} 
          C_{i1} &  C_{i2} &  O_4    & \cdots  &  O_4    \\ 
          C_{i2} &  C_{i1} &  C_{i2} & \ddots  & \vdots  \\ 
          O_4    &  C_{i2} & \ddots  & \ddots  &  O_4    \\ 
         \vdots  & \ddots  & \ddots  & \ddots  &  C_{i2} \\ 
          O_4    & \cdots  &  O_4    &  C_{i2} &  C_{i1}   
           \end{array} 
  \right] \ (i=1, 3),$$ 
where 
$$C_{1k} 
:=\left[ \begin{array}{cccc} 
           0   & -da_k & 0    & -db_k \\ 
          da_k &   0   & db_k &   0   \\ 
           0   & -db_k &  0   & -da_k \\ 
          db_k &   0   & da_k &   0 
           \end{array}   
  \right] , \quad 
  C_{3k} 
:=\left[ \begin{array}{cccc} 
           0   &  da_k & 0    &  db_k \\ 
          da_k &   0   & db_k &   0   \\ 
           0   &  db_k &  0   &  da_k \\ 
          db_k &   0   & da_k &   0 
           \end{array}   
  \right]$$ 
and $a_k$, $b_k$ are functions on $\mbox{\boldmath{$R$}}^m$ ($k=1, 2$). 
Then $D_{11}$, $D_{31}$ satisfy 
${}^t\!D_{11} =-D_{11}$, ${}^t\!D_{31} =D_{31}$, 
and noticing 
$$C_{ik} \wedge C_{jk} =O_4 , \quad 
  C_{i1} \wedge C_{j2} +C_{i2} \wedge C_{j1} =O_4 \quad 
    (i, j=1, 3, k=1, 2),$$  
we obtain $D_{i1} \wedge D_{j1} =O_n$ ($i, j=1, 3$). 
Therefore $D_{11}$, $D_{31}$ satisfy \eqref{dD}. 
\end{ex} 

\begin{rem} 
For any of the above examples, 
$\omega$ is exact and represented as $\omega =dx$  
for a function $x$ on $\mbox{\boldmath{$R$}}^m$ 
valued in the Lie algebra of $SO(2n, 2n)$. 
Therefore a suitable frame field $e$ is given by $\overline{e} \exp (x)$ 
for an ordered pseudo-orthonormal parallel frame field $\overline{e}$ 
of $E$ 
giving the orientation of $E$. 
\end{rem} 

If we suppose $\mu =-$, 
then we have similar discussions and examples. 
In addition, we can have similar discussions and examples 
of $-$-paracomplex structures. 

\begin{rem}\label{rem:rem} 
Suppose $m=4n$ 
and that $E$ is the tangent bundle $T\!E^{4n}_{2n}$ of $E^{4n}_{2n}$. 
Let $h$ be the metric of $E^{4n}_{2n}$ 
and $\nabla$ the Levi-Civita connection of $h$. 
Then the curvature tensor $R$ of $\nabla$ vanishes. 
Therefore we can find 
examples of almost $\varepsilon$-paracomplex structures  
of $E^{4n}_{2n}$ such that the covariant derivatives are nowhere zero 
and satisfy {\rm (a)} in Theorem~\ref{thm:hatnabla0g}. 
\end{rem} 

\section{Frame fields of neutral vector bundles 
of rank 4 with flat metric connections}\label{sect:flatdim4}  

\setcounter{equation}{0} 

In the following, suppose $n=1$. 
Let $e=(e_1 , e_2 , e_3 , e_4 )$ be 
an ordered pseudo-orthonormal frame field of $E$ 
giving the orientation of $E$ such that 
$\hat{\nabla} \Omega_{\varepsilon , 2}$ ($\varepsilon \in \{ +, -\}$) is 
fully light-like. 
Let $\omega =[\omega^i_j ]$ be the connection form of $\nabla$ 
with respect to $e$. 
Then we have \eqref{omegaij}. 
Since $R$ vanishes, $d\omega +\omega \wedge \omega$ vanishes. 
Since $\hat{\nabla} \Omega_{\varepsilon , 2}$ is fully light-like, 
we have \eqref{32+41} for $\mu =\mu_{\varepsilon} \in \{ +, -\}$. 
We set $\psi =[\psi^i_j ] :=\omega \wedge \omega$. 
Then from \eqref{32+41}, we obtain $\psi^4_2 =\varepsilon \psi^3_1$. 
Since $d\omega +\psi$ vanishes, 
we obtain $d(\omega^4_2 -\varepsilon \omega^3_1 )=0$. 
Therefore noticing Proposition~\ref{pro:O-1O+3hor}, 
we can choose $e$ so that $\Omega_0$ 
for $\Omega =\Omega_{\varepsilon , 2}$ is 
horizontal. 
If $e$ is such a frame field, 
then referring to the proof of Proposition~\ref{pro:O-1O+3hor}, 
we obtain $\omega^4_2 -\varepsilon \omega^3_1 =0$. 
However, in the following, 
we suppose only $d(\omega^4_2 -\varepsilon \omega^3_1 )=0$. 
Then there exists a function $f^{\varepsilon}$ 
on $\mbox{\boldmath{$R$}}^m$ 
satisfying $\omega^4_2 -\varepsilon \omega^3_1 =-df^{\varepsilon}$. 
Applying this to $\psi$ and noticing that $d\omega +\psi$ vanishes, 
we obtain 
$$d(\omega^3_2 +\varepsilon \omega^4_1 ) 
   =\mu df^{\varepsilon} \wedge (\omega^3_2 +\varepsilon \omega^4_1 ),$$ 
i.e., 
$d(e^{-\mu f^{\varepsilon}} (\omega^3_2 +\varepsilon \omega^4_1 ))=0$. 
This means that there exists a function $g^{\varepsilon }$ 
on $\mbox{\boldmath{$R$}}^m$ satisfying $dg^{\varepsilon} \not= 0$ and 
$$dg^{\varepsilon} 
=    e^{-\mu f^{\varepsilon}} (\omega^3_2 +\varepsilon \omega^4_1 ) 
=\mu e^{-\mu f^{\varepsilon}} (\omega^4_3 +\varepsilon \omega^2_1 ).$$ 
Therefore we obtain 
\begin{equation} 
\begin{split}
\omega^4_1 & = \varepsilon e^{\mu f^{\varepsilon}} dg^{\varepsilon} 
              -\varepsilon \omega^3_2 , \quad 
\omega^4_2   = \varepsilon \omega^3_1 -df^{\varepsilon} , \\ 
\omega^4_3 & = \mu e^{\mu f^{\varepsilon}} dg^{\varepsilon} 
              -\varepsilon \omega^2_1 . 
\end{split}
\label{omega414243+} 
\end{equation} 
Then the condition that $d\omega +\omega \wedge \omega$ vanishes is 
rewritten into a system of three equations as follows: 
\begin{equation} 
\begin{split} 
d\omega^2_1 = &   2\omega^3_1 \wedge \omega^3_2 
                  -e^{\mu f^{\varepsilon}} 
                   \omega^3_1 \wedge dg^{\varepsilon} \\ 
              &   +\varepsilon \omega^3_2 \wedge df^{\varepsilon} 
                  +\varepsilon e^{\mu f^{\varepsilon }} 
                   df^{\varepsilon} \wedge dg^{\varepsilon} , \\ 
d\omega^3_1 = &  2\omega^2_1 \wedge \omega^3_2 
                 -e^{\mu f^{\varepsilon}} 
                  \omega^2_1 \wedge dg^{\varepsilon} 
                 +\varepsilon \mu e^{\mu f^{\varepsilon}} 
                  \omega^3_2 \wedge dg^{\varepsilon} , \\ 
d\omega^3_2 = & -2\omega^2_1 \wedge \omega^3_1 
                + \varepsilon \omega^2_1 \wedge  df^{\varepsilon} \\ 
              & - \varepsilon \mu e^{\mu f^{\varepsilon}} 
                  \omega^3_1 \wedge dg^{\varepsilon} 
                + \mu e^{\mu f^{\varepsilon}} 
                  df^{\varepsilon} \wedge dg^{\varepsilon} . 
\end{split} 
\label{omega213132+} 
\end{equation} 

Let $f^{\varepsilon}$, $g^{\varepsilon}$ be functions 
on $\mbox{\boldmath{$R$}}^m$ 
with $dg^{\varepsilon} \not= 0$. 
Let $\omega^2_1$, $\omega^3_1$, $\omega^3_2$ be $1$-forms 
on $\mbox{\boldmath{$R$}}^m$ satisfying \eqref{omega213132+}. 
Let $\omega^4_1$, $\omega^4_2$, $\omega^4_3$ be as in \eqref{omega414243+} 
and $\omega$ a $4\times 4$ matrix such that the $(i, j)$-component is 
given by $\omega^i_j$ with \eqref{omegaij}. 
Then $d\omega +\omega \wedge \omega$ vanishes. 
Therefore, noticing $R=0$, 
we see that there exists 
an ordered pseudo-orthonormal frame field $e=(e_1 , e_2 , e_3 , e_4 )$ 
of $E$ giving the orientation of $E$ and satisfying $\nabla e=e\omega$. 
Then for $\Omega_{\varepsilon , 2}$ as in \eqref{Omega+-123}, 
we obtain $\hat{\nabla} \Omega_{\varepsilon , 2} 
          =\alpha \otimes \Omega_0$, 
where 
\begin{equation} 
\alpha   :=e^{\mu f^{\varepsilon}} dg^{\varepsilon} , \quad 
\Omega_0 :=\Omega_{-\varepsilon , 1} +\mu \Omega_{\varepsilon  , 3} . 
\label{aOE42} 
\end{equation} 
Since $dg^{\varepsilon} \not= 0$, we have $\alpha \not= 0$. 
Since $d(\omega^4_2 -\varepsilon \omega^3_1 )=0$, 
we see from Proposition~\ref{pro:O-1O+3hor} 
that $\Omega_0$ can be horizontal for a suitable frame field. 
Let $J_{\varepsilon}$ be an $\varepsilon$-paracomplex structure of $E$ 
corresponding to $\Omega_{\varepsilon , 2}$. 
Then we have $J_{\varepsilon} e_1 =             e_3$, 
             $J_{\varepsilon} e_2 =-\varepsilon e_4$, and 
$J=J_{\varepsilon}$ satisfies \eqref{nablaJ} 
with \eqref{xitheta} and \eqref{aOE42}. 

Hence we obtain 

\begin{pro}\label{pro:gsec} 
The following two conditions are equivalent$:$ 
\begin{itemize} 
\item[{\rm (a)}]{the covariant derivative of $\Omega_{\varepsilon , 2}$ is 
fully light-like\/$;$} 
\item[{\rm (b)}]{the connection form $\omega =[\omega^i_j ]$ of $\nabla$ 
with respect to the frame field $e=(e_1 , e_2 , e_3 , e_4 )$ 
satisfies \eqref{omega414243+}, \eqref{omega213132+} 
and $dg^{\varepsilon} \not= 0$ 
for functions $f^{\varepsilon}$, $g^{\varepsilon}$ 
on $\mbox{\boldmath{$R$}}^m$ and $\mu \in \{ +, -\}$.} 
\end{itemize} 
In addition, if one of {\rm (a)}, {\rm (b)} holds, 
then $\hat{\nabla} \Omega_{\varepsilon , 2}$ is represented 
as $\hat{\nabla} \Omega_{\varepsilon , 2} =\alpha \otimes \Omega_0$, 
where $\alpha$, $\Omega_0$ are as in \eqref{aOE42} 
and $\Omega_0$ can be horizontal for a suitable frame field. 
\end{pro} 

Suppose $m=4$ 
and that $E$ is the tangent bundle $T\!E^4_2$ of $E^4_2$. 
Let $h$ be the metric of $E^4_2$ 
and $\nabla$ the Levi-Civita connection of $h$. 
Then the curvature tensor $R$ of $\nabla$ vanishes. 
Therefore, referring to the above discussions, we obtain 

\begin{cor}\label{cor:E42J} 
For $E=T\!E^4_2$, {\rm (a)}, {\rm (b)} in Proposition~\ref{pro:gsec} are 
equivalent to each other. 
In addition, 
if one of {\rm (a)}, {\rm (b)} in Proposition~\ref{pro:gsec} holds, 
then the remaining statements in Proposition~\ref{pro:gsec} hold. 
\end{cor} 

Suppose that there exist sections $\Omega_{\pm}$ 
of $U_-\!\left(\bigwedge^2_{\pm}\!E\right)$ respectively 
such that the covariant derivatives are fully light-like. 
Then noticing 
the double covering $SO_0(2,2)\longrightarrow SO_0 (1,2)\times SO_0(1,2)$, 
we see that there exists 
an ordered pseudo-orthonormal frame field $e=(e_1 , e_2 , e_3 , e_4 )$ 
of $E$ 
such that both of $\hat{\nabla} \Omega_{\pm , 2}$ are fully light-like. 
Let $\omega =[\omega^i_j ]$ be the connection form of $\nabla$ 
with respect to $e$. 
Then we have \eqref{32+41} for $\mu =\mu_{\varepsilon}$, 
$\varepsilon =+$, $-$. 
Therefore we have the second relation in \eqref{omega414243+} 
for $\varepsilon =+$, $-$. 
Therefore we have 
\begin{equation} 
\omega^3_1 = \dfrac{1}{2} (df^+ -df^-), \quad 
\omega^4_2 =-\dfrac{1}{2} (df^+ +df^-). 
\label{3142+-} 
\end{equation} 
Noticing \eqref{32+41}, 
we obtain either $\omega^4_1 =\mu \omega^2_1$, 
                 $\omega^4_3 =\mu \omega^3_2$ 
              or $\omega^3_2 =\mu \omega^2_1$, 
                 $\omega^4_3 =\mu \omega^4_1$ for $\mu \in \{ +, -\}$. 
Suppose $\omega^4_1 =\mu \omega^2_1$ 
    and $\omega^4_3 =\mu \omega^3_2$.  
Then by \eqref{omega414243+}, we obtain 
\begin{equation} 
\begin{split} 
 \omega^2_1 & 
=\dfrac{\mu}{2} 
 \left( e^{\mu f^+} dg^+ -e^{\mu f^-} dg^- 
 \right) , \\  
 \omega^3_2 & 
=\dfrac{1}{2} 
 \left( e^{\mu f^+} dg^+ +e^{\mu f^-} dg^- 
 \right) 
\end{split} 
\label{dfdg213132+} 
\end{equation} 
and 
\begin{equation} 
\begin{split} 
 \omega^4_1 & 
=\dfrac{1}{2} 
 \left( e^{\mu f^+} dg^+ -e^{\mu f^-} dg^- 
 \right) , \\ 
 \omega^4_3 & 
=\dfrac{\mu}{2} 
 \left( e^{\mu f^+} dg^+ +e^{\mu f^-} dg^- 
 \right) . 
\end{split} 
\label{dfdg414243+} 
\end{equation} 
Suppose $\omega^3_2 =\mu \omega^2_1$ 
    and $\omega^4_3 =\mu \omega^4_1$. 
Then by \eqref{omega414243+}, we obtain 
\begin{equation} 
\begin{split} 
 \omega^2_1 & 
=\dfrac{\mu}{2} 
 \left( e^{\mu f^+} dg^+ +e^{-\mu f^-} dg^- 
 \right) , \\ 
 \omega^3_2 & 
=\dfrac{1}{2} 
 \left( e^{\mu f^+} dg^+ +e^{-\mu f^-} dg^- 
 \right) 
\end{split} 
\label{dfdg213132-} 
\end{equation} 
and 
\begin{equation} 
\begin{split} 
 \omega^4_1 & 
=\dfrac{1}{2} 
 \left( e^{\mu f^+} dg^+ -e^{-\mu f^-} dg^- 
 \right) , \\ 
 \omega^4_3 & 
=\dfrac{\mu}{2} 
 \left( e^{\mu f^+} dg^+ -e^{-\mu f^-} dg^- 
 \right) . 
\end{split} 
\label{dfdg414243-} 
\end{equation} 

Let $f^{\pm}$, $g^{\pm}$ be functions on $\mbox{\boldmath{$R$}}^m$ 
with $dg^{\pm} \not= 0$ 
and $\omega^3_1$, $\omega^4_2$ 1-forms on $\mbox{\boldmath{$R$}}^m$ 
given by \eqref{3142+-}. 
Let $\omega^2_1$, $\omega^3_2$, $\omega^4_1$, $\omega^4_3$ be 
1-forms on $\mbox{\boldmath{$R$}}^m$ given by 
either \eqref{dfdg213132+}, \eqref{dfdg414243+} 
or     \eqref{dfdg213132-}, \eqref{dfdg414243-}. 
Then $\omega^2_1$, $\omega^3_1$, $\omega^3_2$, 
     $\omega^4_1$, $\omega^4_2$, $\omega^4_3$ satisfy 
\eqref{omega414243+}, \eqref{omega213132+} 
for $\varepsilon =+$, $-$ and $\mu =\mu_{\varepsilon}$. 
Let $\omega$ be a $4\times 4$ matrix such that the $(i, j)$-component is 
given by $\omega^i_j$ with \eqref{omegaij}. 
Then there exists 
an ordered pseudo-orthonormal frame field $e=(e_1 , e_2 , e_3 , e_4 )$ 
of $E$ with $\nabla e=e\omega$ giving the orientation of $E$ and 
both of $\hat{\nabla} \Omega_{\pm , 2}$ are fully light-like. 

Hence we obtain 

\begin{thm}\label{thm:gsec+-} 
The following two conditions are equivalent$:$ 
\begin{itemize} 
\item[{\rm (a)}]{the covariant derivatives of both of $\Omega_{\pm , 2}$ 
are fully light-like\/$;$} 
\item[{\rm (b)}]{the connection form $\omega =[\omega^i_j ]$ of $\nabla$ 
with respect to the frame field $e=(e_1 , e_2 , e_3 , e_4 )$ 
satisfies \eqref{3142+-} and 
either \eqref{dfdg213132+}, \eqref{dfdg414243+} 
or     \eqref{dfdg213132-}, \eqref{dfdg414243-} 
for functions $f^{\pm}$, $g^{\pm}$ on $\mbox{\boldmath{$R$}}^m$ 
with $dg^{\pm} \not= 0$.} 
\end{itemize} 
\end{thm} 

\begin{rem}\label{rem:rem2_0}
Based on Theorem~\ref{thm:gsec+-}, 
we see that each pair of almost $\pm$-paracomplex structures 
of $E^4_2$ giving sections $\Omega_{\pm , 2}$ 
of the two time-like twistor spaces as in (a) in Theorem~\ref{thm:gsec+-} 
is given by functions $f^{\pm}$, $g^{\pm}$ 
on $E^4_2$ with $dg^{\pm} \not= 0$. 
\end{rem} 

\begin{rem}\label{rem:rem2} 
Suppose that the covariant derivatives of $\Omega_{\pm , 2}$ are 
fully light-like.  
Then as was already seen, we have 
either $\omega^4_1 =\mu \omega^2_1$, 
       $\omega^4_3 =\mu \omega^3_2$ 
    or $\omega^3_2 =\mu \omega^2_1$, 
       $\omega^4_3 =\mu \omega^4_1$ for $\mu \in \{ +, -\}$. 
The former (respectively, latter) case corresponds to 
(i) (respectively, (ii)) in the last paragraph 
of Section~\ref{sect:cgmtms}. 
\end{rem} 

\section*{Acknowledgements} 

The author is grateful to Professor Johann Davidov 
for helpful communications. 
This work was supported by 
JSPS KAKENHI Grant Number JP21K03228.

\vspace{4mm} 

\par\noindent 
\footnotesize{Faculty of Advanced Science and Technology, 
              Kumamoto University \\ 
              2--39--1 Kurokami, Chuo-ku, Kumamoto 860--8555 Japan} 

\par\noindent  
\footnotesize{E-mail address: andonaoya@kumamoto-u.ac.jp} 


\begin{thebibliography}{99} 

\bibitem{AGK}D.\,V.\,Alekseevsky, B.\,Guilfoyle and W.\,Klingenberg, 
On the geometry of spaces of oriented geodesics, 
Ann. Global Anal. Geom. {\bf 40} (2011) 389--409. 
\textit{Erratum\/}: Ann. Global Anal. Geom. {\bf 50} (2016) 97--99. 

\bibitem{ando4}N.\,Ando, 
Surfaces with zero mean curvature vector in neutral $4$-manifolds, 
Diff. Geom. Appl. {\bf 72} (2020) 101647, 31 pp. 

\bibitem{ando7}N.\,Ando, 
Nilpotent structures of neutral 4-manifolds and light-like surfaces, 
Developments in Lorentzian Geometry, 
Springer Proceedings in Mathematics \& Statistics {\bf 389}, 
Springer, 2022, 
13--28. 

\bibitem{ando5}N.\,Ando, 
The lifts of surfaces in neutral 4-manifolds into the 2-Grassmann bundles, 
Diff. Geom. Appl. {\bf 91} (2023) 102073, 25 pp. 

\bibitem{ando6}N.\,Ando, 
Time-like surfaces with zero mean curvature vector 
in four-dimensional neutral space forms, 
Proc. Int. Geom. Cent. {\bf 17} (2024) 36--55. 

\bibitem{ando8}N.\,Ando, 
Nilpotent structures of oriented neutral vector bundles, 
preprint; arXiv:2405.05003. 

\bibitem{AHHK}
N.\,Ando, K.\,Hamada, K.\,Hashimoto and S.\,Kato, 
Regularity of ends of zero mean curvature surfaces 
in $\mathbf{R}^{2,1}$, 
J. Math. Soc. Japan {\bf 74} (2022) 1295--1334. 

\bibitem{AK}N.\,Ando and T.\,Kihara, 
Horizontality in the twistor spaces associated with vector bundles of 
rank 4 on tori, J. Geom. {\bf 112} (2021) 19, 26 pp. 

\bibitem{asgeirsson}L.\,\'{A}sgeirsson, 
\"{U}ber eine Mittelwertseigenschaft von L\"{o}sungen homogener 
linearer partieller Differentialgleichungen 2. 
Ordnung mit konstanten Koeffizienten, 
Math. Ann. {\bf 113} (1937) 321--346. 

\bibitem{BDM}D.\,Blair, J.\,Davidov and O.\,Mu\v{s}karov, 
Isotropic K\"{a}hler hyperbolic twistor spaces, 
J. Geom. Phys. {\bf 52} (2004) 74--88. 

\bibitem{BDM2}D.\,Blair, J.\,Davidov and O.\,Mu\v{s}karov, 
Hyperbolic twistor spaces, Rocky Mountain J. Math. {\bf 35} (2005) 
1437--1465. 

\bibitem{BGGNV}M.\,Brozos-V\'azquez, 
E.\,Garc\'{\i}a-R\'{\i}o, 
P.\,Gilkey, 
S.\,Nik\v{c}evi\'c, 
R.\,V\'{a}zquez-Lorenzo, 
The geometry of Walker manifolds, 
Synthesis Lectures on Mathematics and Statistics, 
Morgan and Claypool, 2009. 

\bibitem{CG}G.\,Cobos and B.\,Guilfoyle, 
An extension of Asgeirsson's mean value theorem 
for solutions of the ultra-hyperbolic equation in dimension four, 
Diff. Geom. Appl. {\bf 79} (2021) 101795, 17 pp. 

\bibitem{DDGMMV}J.\,Davidov, 
J.\,C.\,D\'{\i}az-Ramos, 
E.\,Garc\'{\i}a-R\'{\i}o, 
Y.\,Matsushita, 
O.\,Mu\v{s}karov and 
R.\,V\'{a}zquez-Lorenzo, 
Almost K\"{a}hler Walker 4-manifolds, 
J. Geom. Phys. {\bf 57} (2007) 1075--1088. 

\bibitem{DDGMMV2}J.\,Davidov, 
J.\,C.\,D\'{\i}az-Ramos, 
E.\,Garc\'{\i}a-R\'{\i}o, 
Y.\,Matsushita, 
O.\,Mu\v{s}karov and
R.\,V\'{a}zquez-Lorenzo, 
Hermitian-Walker 4-manifolds, 
J. Geom. Phys. {\bf 58} (2008) 307--323. 

\bibitem{dunajski0}M.\,Dunajski, 
Anti-self-dual four-manifolds with a parallel real spinor, 
R. Soc. Lond. Proc. Ser. A Math. Phys. Eng. Sci. {\bf 458} (2002) 1205--1222. 

\bibitem{dunajski}M.\,Dunajski, 
Null K\"{a}hler geometry and isomonodromic deformations, 
Commun. Math. Phys. {\bf 391} (2022) 77--105. 

\bibitem{ES}J.\,Eells and S.\,Salamon, 
Twistorial construction of harmonic maps of surfaces 
into four-manifolds, 
Annali della Scuola Normale Superiore di Pisa, Classe di Scienze 
{\bf 12} (1985) 589--640. 

\bibitem{friedrich}T.\,Friedrich, 
On surfaces in four-spaces, 
Ann. Glob. Anal. Geom. {\bf 2} (1984) 257--287. 

\bibitem{GM}E.\,Garc\'{\i}a-R\'{\i}o and Y.\,Matsushita, 
Isotropic K\"{a}hler structures on Engel 4-manifolds, 
J. Geom. Phys. {\bf 33} (2000) 288--294. 

\bibitem{GG}N.\,Georgiou and B.\,Guilfoyle, 
Almost paracomplex structures on 4-manifolds, 
Diff. Geom. Appl. {\bf 82} (2022) 101890, 24 pp. 

\bibitem{guilfoyle}B.\,Guilfoyle, 
From CT scans to 4-manifold topology via neutral geometry, 
Irish Math. Soc. Bull. {\bf 91} (2023) 9--32. 

\bibitem{GK}B.\,Guilfoyle and W.\,Klingenberg, 
An indefinite K\"{a}hler metric on the space of oriented lines, 
J. London Math. Soc. {\bf 72} (2005) 497--509. 

\bibitem{HM}K.\,Hasegawa and K.\,Miura, 
Extremal Lorentzian surfaces with null $\tau$-planar geodesics 
in space forms, 
Tohoku Math. J. {\bf 67} (2015) 611--634. 

\bibitem{JR}G.\,Jensen and M.\,Rigoli, 
Neutral surfaces in neutral four-spaces, 
Matematiche (Catania) {\bf 45} (1990) 407--443. 

\bibitem{john}F.\,John, 
The ultrahyperbolic differential equation with four independent variables, 
Duke Math. J. {\bf 4} (1938) 300--322. 

\bibitem{pavsic}M.\,Pav\v{s}i\v{c}, 
Quantum field theories in spaces with neutral signatures, 
J. Phys.: Conf. Ser. {\bf 437} (2013) 012006, 29 pp. 

\bibitem{salvai}M.\,Salvai, 
On the geometry of the space of oriented lines of the hyperbolic space, 
Glasg. Math. J. {\bf 49} (2007) 357--366. 

\bibitem{salvai2}M. Salvai, 
Global smooth fibrations of $\mathbb{R}^3$ by oriented lines, 
Bull. Lond. Math. Soc. {\bf 41} (2009) 155--163. 

\bibitem{walker}A.\,G.\,Walker, 
Canonical form for a Riemannian space with a parallel field of null 
planes, Quart. J. Math. Oxford (2) {\bf 1} (1950) 69--79. 

\end{thebibliography}
\end{document}